\documentclass[12pt]{amsart}

\usepackage{amsmath,amssymb,amsfonts,amsthm}
\usepackage{commath}
\usepackage{epsfig}
\usepackage{lipsum}
\usepackage{graphicx}
\usepackage{calc}
\usepackage{gensymb}
\usepackage{wrapfig}
\usepackage[sort,numbers]{natbib}
\usepackage{placeins}
\usepackage[dvipsnames]{xcolor}
\usepackage[normalem]{ulem}
\usepackage{bbm}
\usepackage{adjustbox}
\usepackage{float}
\usepackage{xkeyval}
\usepackage{mathtools}
\usepackage{tikz}
\usepackage{subcaption}
\usepackage{hyperref}

\usetikzlibrary{
    shapes.geometric,
    trees,
    arrows,
    calc,
    positioning,
    decorations.markings,
    arrows.meta
}

\tikzset{
    oriented/.style={
        decoration={
            markings,
            mark=at position 0.5 with {\arrow{>}}
        },
        postaction={decorate}
    },
    strand/.style={
        line width=.8pt,
        line cap=round
    },
    band/.style={
        draw=gray!55,
        line width=2mm,
        line cap=round
    },
    >={Stealth[length=2mm]}
}

\renewcommand{\thefootnote}{}

\theoremstyle{plain}
\newtheorem{theorem}{\indent\sc Theorem}[section]
\newtheorem{lemma}[theorem]{\indent\sc Lemma}
\newtheorem{corollary}[theorem]{\indent\sc Corollary}
\newtheorem{proposition}[theorem]{\indent\sc Proposition}

\theoremstyle{definition}
\newtheorem{definition}[theorem]{\indent\sc Definition}
\newtheorem{remark}[theorem]{\indent\sc Remark}

\newtheorem{example}[theorem]{\indent\sc Example}

\newtheorem*{notation}{\indent\sc Notation}

\newcommand{\D}{\mathbbm{D}}

\DeclareMathOperator{\mdet}{mdet}

\title{Determinant and Quasi-Alternateness of Checkerboard Colorable Links}

\author{Mohammed Sabak}
\address{Fundamental and Applied Mathematics Laboratory, Hassan II University EST, Casablanca, Morocco}
\email{mohammed.sabak-etu@etu.univh2c.ma}

\author{Khaled Qazaqzeh}
\address{Department of Mathematics, Faculty of Science, Yarmouk University, Irbed, Jordan, 21163}
\email{qazaqzeh@yu.edu.jo}

\author{Hamid Abchir}
\address{Fundamental and Applied Mathematics Laboratory, Hassan II University, Casablanca, Morocco}
\email{hamid.abchir@univh2c.ma}

\date{}

\begin{document}

\maketitle

\begingroup
\renewcommand{\thefootnote}{}
\endgroup

\begin{abstract}
We give a connected state expansion of the determinant of checkerboard colorable links and extend the notion of quasi-alternateness from classical links to this setting. We provide an obstruction to quasi-alternateness and we characterize classicality of alternating virtual links. Then, we show that quasi-alternateness is preserved under rational extensions of quasi-alternating crossings. Finally, we define the family of almost alternating virtual links. We show that extending the dealternator of an almost alternating virtual diagram of a quasi-alternating checkerboard colorable link preserves the quasi-alternateness.
\end{abstract}
\section{Introduction}
The theory of virtual knots was introduced by Kauffman \cite{K1,K}, motivated by the study of knots in thickened surfaces and abstract Gauss codes. Goussarov, Polyak, and Viro \cite{GPV} proved that the natural map from classical links to virtual links is injective. Consequently, virtual knot theory provides a natural extension of classical knot theory. Therefore, it is natural to seek extensions of classical knot invariants to the virtual setting. While many invariants admit natural extensions, others either do not extend naturally or admit several inequivalent extensions, as discussed in \cite[Section~1.5]{GPV}. The determinant is one such invariant, and several extensions have been introduced in \cite{imlee2010,boden2022,boden2023,abchir2026}. In this paper, for a checkerboard colorable link $L$, we consider the determinant introduced by Boden and Karimi \cite{boden2022}  denoted by $\det(L)$. This invariant extends the classical determinant to the class of checkerboard colorable links.

For classical links, the determinant admits several equivalent formulations, including those in terms of the Goeritz matrix, the Alexander polynomial, and the Jones polynomial. In particular, the absolute value of the determinant of a Goeritz matrix gives the determinant of a classical link. The Goeritz matrix construction can be adapted to checkerboard colorable links \cite{imlee2010}. In the virtual setting, this construction gives a pair of Goeritz matrices and corresponding determinants, which depend on the choice of spanning surface and, in general, do not coincide with the determinant of a checkerboard colorable link \cite{boden2022-2}. The mock Seifert matrix introduced by Boden and Karimi in \cite{boden2023} instead provides a formulation of the determinant appropriate to the virtual setting. The Alexander and Jones polynomial formulations are given by the evaluations at $-1$: If $\Delta_L(t)$ and $J_L(t)$ denote the Alexander and Jones polynomials of a classical link $L$, respectively, then $\det(L)=|\Delta_L(-1)|=|J_L(-1)|$. This identity, together with the spanning-tree expansion of the Jones polynomial \cite{thistlethwaite1987}, was used by Champanerkar and Kofman \cite{kofman2009} to derive the following spanning-tree expansion of the determinant:
\begin{align}
\det(L)=\left|\sum_v(-1)^v s_v(D)\right|,
\label{formula}
\end{align}
where $s_v(D)$ denotes the number of spanning trees of the Tait graph of a diagram $D$ having $v$ positive edges. This combinatorial formulation has been used in various aspects of the study of quasi-alternating links, notably to suggest bounds on the crossing number \cite{qazaqzeh2013}, to show that the quasi-alternating property is preserved under tangle extensions of certain crossings \cite{kofman2009,abchir2020,chbili2020,kaur2026II,kaur2026III}, and to contribute to the characterization of quasi-alternating links within certain families \cite{kofman2009,qazaqzeh2015,champanerkar2015,issa2018}.

In contrast, the corresponding situation for checkerboard colorable links is substantially different. The Alexander polynomial is defined only for the subclass of almost classical links \cite{boden2017} and the Jones-Kauffman polynomial does not yield the determinant in general as observed in \cite{krb,abchir2026}.

The main purpose of this paper is to establish a combinatorial formula for the determinant of checkerboard colorable links analoguous to Equation \eqref{formula} that appears in the proof of \cite[Theorem\,1]{kofman2009}. More precisely, we prove the following:

\begin{theorem}
Let $L$ be a virtual link admitting a checkerboard colorable diagram $D$. Then $$
\det(L)=\abs{\displaystyle\sum_{s\in\mathcal{S}_D}(-1)^{a_s}},$$
where $\mathcal{S}_D$ is the set of connected states of $D$, and $a_s$ denotes the number of markers labelled $a$ in the connected state $s$.
\label{th1}
\end{theorem}

Although the resulting formula resembles its classical counterpart (Equation \ref{formula}), its proof cannot be obtained by adapting the classical argument, which relies on the relationship between the determinant and the Jones polynomial in the case of classical links that does not extend to checkerboard colorable links.

We use this formula to generalize the determinant formulas for several families of classical links to their checkerboard colorable analogues  (Propositions \ref{prop3.11}, \ref{prop2}, and \ref{prop3.10}). We also compare the determinant with the mock-determinant (Proposition \ref{prop4.1}), another extension of the classical determinant to the virtual setting introduced in \cite{abchir2026}. The mock-determinant of a virtual link $L$, denoted by $\mdet(L)$, is defined as the $L_1$-norm of the Jones-Kauffman polynomial evaluated at the eighth root of unity $e^{\frac{i\pi}{4}}$. This comparison yields an obstruction to alternateness of virtual links and a characterization of classicality for alternating virtual links.  More precisely, we prove the following theorem:

\begin{theorem}
 If $L$ is an alternating virtual link, then $L$ is classical if and only if $\mdet(L) = \det(L)$.
\label{th2}
\end{theorem}

Finally, we extend the notion of quasi-alternateness to checkerboard colorable virtual links and note that every alternating virtual link is quasi-alternating. Furthermore, we give an obstruction to quasi-alternateness. In particular, we prove the following theorem:
\begin{theorem}
If $L$ is a quasi-alternating checkerboard colorable link, then $\mdet(L) \leq \det(L)$. 
\label{th3}
\end{theorem}

Using our combinatorial determinant formula, we show that the twisting technique introduced in \cite{kofman2009} preserves quasi-alternateness and produces an infinite family of quasi-alternating checkerboard colorable links when applied at any crossing where the link is quasi-alternating (Theorem \ref{new}). By the same way, we prove that this technique provides an infinite family of quasi-alternating checkerboard colorable links when applied at the dealternator of an almost alternating link diagram representing a quasi-alternating link (Theorem \ref{th5.5}). 

The paper is organized as follows. Section \ref{sec2} recalls the preliminaries needed for the proofs. In Section \ref{sec3}, we prove Theorem \ref{th1} and apply the resulting formula to generalize determinant formulas for several families of classical links to their checkerboard colorable analogues. Section \ref{sec4} is devoted to the comparison between the determinant and the mock-determinant and to the proof of Theorem \ref{th2}. In Section \ref{sec5}, we extend quasi-alternateness to checkerboard colorable links and prove Theorem \ref{th3}.

\section{Preliminaries}
\label{sec2}
In this section, we recall the basic concepts and terminologies required to state and prove the main results of this paper.

\subsection{Virtual links}
\begin{figure}
\centering
\includegraphics[scale=0.3]{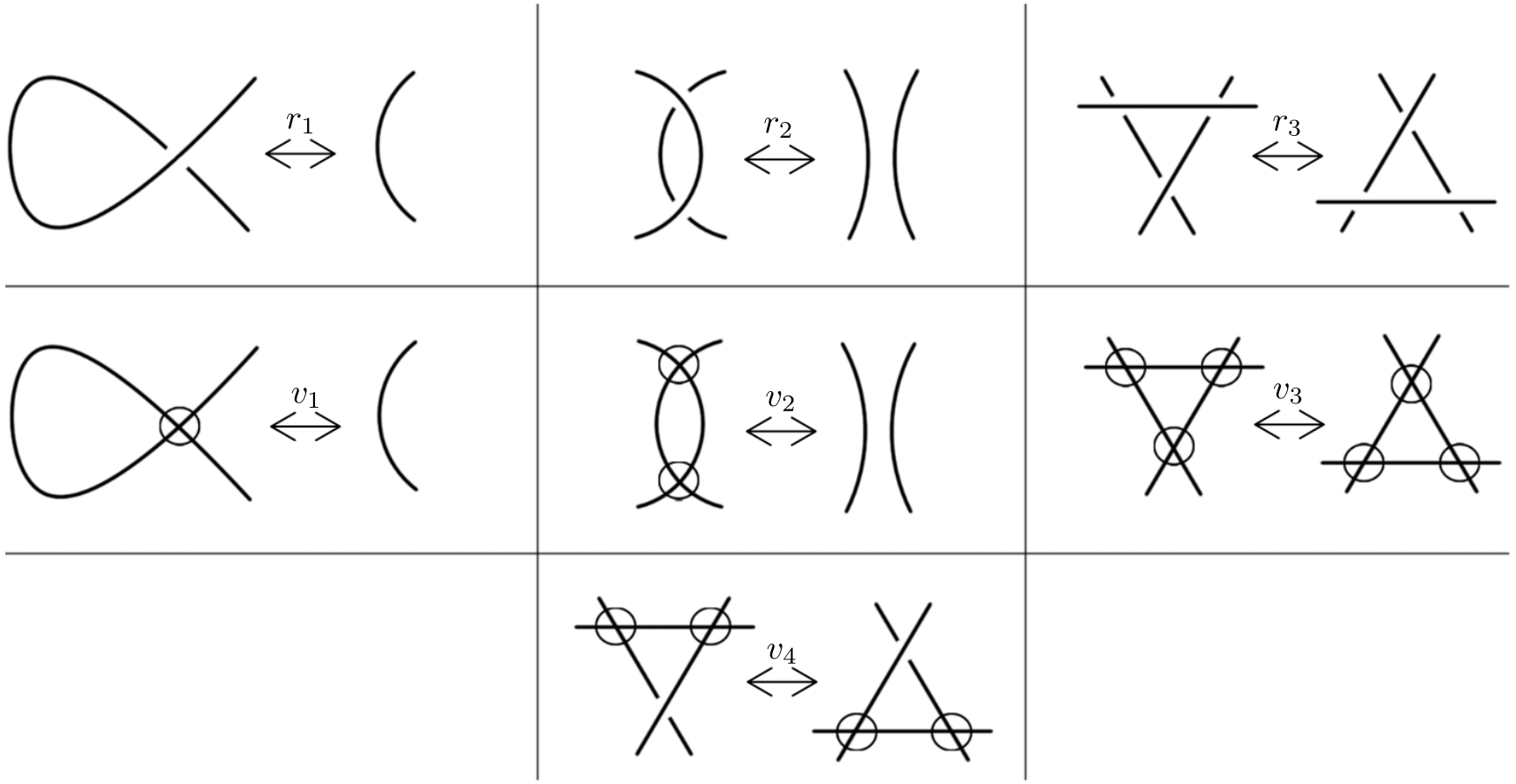}
\caption{}
\label{grm}
\end{figure}

We recall the definition of a virtual link diagram and a virtual link. A \textit{virtual link diagram} is a generic immersion of one or several circles in the plane with only double points, such that each double point is either a \textit{classical crossing} indicated by over and under-crossings, or a \textit{virtual crossing} indicated by a circle. Two virtual link diagrams are \textit{equivalent} if they are related by a finite sequence of the generalized Reidemeister moves which are local modifications of virtual link diagrams depicted in Figure \ref{grm}. A \textit{virtual link} is an equivalence class of virtual diagrams under the generalized Reidemeister moves. Throughout the paper, the term “crossing” refers to a classical crossing unless otherwise specified.

\subsection{Abstract link diagrams and checkerbaord colorability of links}\label{section2.1}

We recall the notion of \textit{abstract link diagram} that was first introduced by N. Kamada and S. Kamada in \cite{kk}. An abstract link diagram is a pair $P=(S,d)$ consisting of a compact oriented surface $S$ and a link diagram $d$ on $S$ that is a deformation retract of $S$. If there exists an orientation preserving embedding $f : S \rightarrow F$ into a closed connected oriented surface $F$, we say that $F$ is a \textit{surface realization} of the abstract link diagram $P=(S,d)$, and that the link diagram $f(d)$ on $F$ is a \textit{link diagram realization} of $P$. 

An abstract link diagram $(S_D,\widetilde D)$ can be associated to any given virtual link diagram $D$ in the plane according to the construction in \cite[Section\,4]{kk}. 
This association induces a bijection between virtual link diagrams modulo generalized Reidemeister moves and abstract link diagrams modulo abstract equivalence as shown in \cite[Thereom\,4.1]{kk}. Moreover, as explained in \cite[Section\,2]{kk}, attaching disks to $S_D$ along its boundary components yields the closed orientable surface $\Sigma_D$ called the Carter surface. In fact, $\Sigma_D$ has the minimum genus among all surface realizations of $(S_D,\widetilde D)$ that is called the \textit{supporting genus} of $(S_D,\widetilde D)$ and is denoted by $sg(D)$. 

A virtual link diagram $D$ is a \textit{split diagram} if it can be transformed into a disconnected diagram using moves $(v_1)-(v_4)$ (Figure \ref{grm}). If $D$ is disconnected, then it can be easily seen that the lift $\mathcal{D}$ of $D$ to the Carter surface $\Sigma_D$ will be also disconnected. In such a case the link $L_D$ in $\Sigma_D \times \left[ 0,1 \right]$ represented by $\mathcal{D}$ is called a \textit{split link}. 

An abstract link diagram $P=(S,d)$ is \textit{checkerboard colorable} if to each region of the complement of the underlying graph $|d|$ in $S$, one can assign one of two colors, say black and white, such that any two adjacent regions sharing an edge of $|d|$ have different colors. A virtual link diagram $D$ is said to be \textit{checkerboard colorable} if the associated abstract link diagram $(S_D,\widetilde D)$ admits a checkerboard coloring. In such a case, the link diagram $\mathcal{D}$ in $\Sigma_D$ and the associated link $L_D$ in $\Sigma_D \times \left[ 0,1 \right]$ are also said to be \textit{checkerboard colorable}. This is equivalent to say that the $1$-cycle represented by the link $L_D$ is null-homologous modulo $2$, and thus the link bounds a (possibly non-orientable) surface in $\Sigma_D \times \left[ 0,1 \right]$. A virtual link is \textit{checkerboard colorable} if it admits a checkerboard colorable diagram. These particular virtual links were introduced by Kamada in \cite{kamada2002}. 

\begin{figure}
\centering
\includegraphics[width=0.6\linewidth]{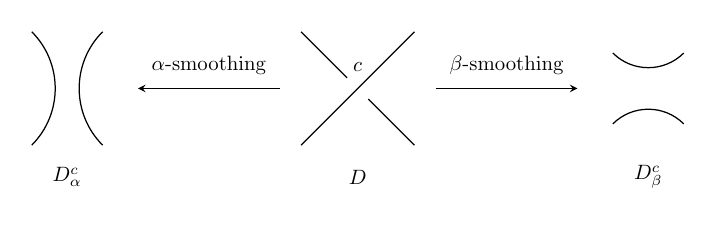}
\caption{The link diagram $D$ and its smoothings $D^{c}_{\alpha}$ and $D^{c}_{\beta}$ at the crossing $c$.}
\label{fig.1}
\end{figure}

A virtual link diagram is \textit{alternable} if it can be made alternating by changing some of its classical crossings. Kamada showed that a virtual link diagram is checkerboard colorable if and only if it is alternable \cite[Lemma\,7]{kamada2002}. In particular, any classical or alternating virtual link diagram is checkerboard colorable. Another characterization of checkerboard colorability of a virtual link diagram $D$ is given by Karimi in his thesis \cite{karimi2018} as follows: Each crossing of $D$ can be resolved into an $\alpha$-smoothing or a $\beta$-smoothing as in Figure \ref{fig.1}. A \textit{state} of $D$ is the virtual link diagram obtained by choosing either an $\alpha$- or a $\beta$-smoothing at each crossing of $D$, leaving only virtual crossings. If $D$ has $n$ crossings, then there are $2^n$ states of $D$. For a given state $s$, the \textit{dual state} $\overline{s}$ is obtained from $s$ by changing all $\alpha$-smoothings to $\beta$-smoothing, and vice versa. Let $(S_D,\widetilde D)$ be the abstract link diagram associated to $D$. We say that $D$ has the \textit{boundary property} if there exists a state $s_{\partial}$ such that $\partial S_D = s_{\partial} \cup \overline{s}_{\partial}$. Karimi showed in \cite[Lemma 3.17]{karimi2018} that a virtual link diagram $D$ is checkerboard colorable if and only if it has the boundary property. 

\begin{figure}
\centering
\includegraphics[width=0.5\linewidth]{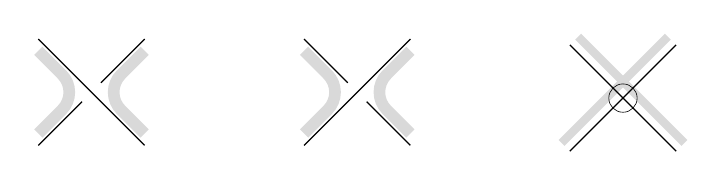}
\caption{}
\label{ccol}
\end{figure}

If $(S_D,\widetilde D)$ is checkeboard colorable, then there is an induced coloring as shown in Figure \ref{ccol} of the sides of \textit{long arcs} of $D$, i.e., strands of $D$ starting and finishing at classical crossings. In such a figure, a band along one side of a regular neighborhood of each long arc is colored so that, at each classical crossing, the coloring switches sides and the colored bands along the over-arc and the under-arc connect within the same region of the complement of the diagram. At a virtual crossing, the coloring passes through without any restriction. See the left of Figure \ref{ccdcs} for an example of a checkerboard colorable virtual link diagram. A notable fact about checkerboard colorability is that it behaves well under generalized Reidemeister moves as shown in \cite[Theorem\,2.9]{im2}. We assume from now on that any virtual link diagram is checkerboard colorable, unoriented and non-split unless mentioned otherwise. 

\subsection{Tangles and rational extensions}
A \textit{virtual tangle diagram} is a generic immersion in the unit disk $D^2$ of a compact $1$-manifold $T$ such that $\partial T$ consists of four distinct points lying on $\partial D^2$ given by 
$
NW=\left(-\frac{1}{\sqrt2},\frac{1}{\sqrt2}\right),
NE=\left(\frac{1}{\sqrt2},\frac{1}{\sqrt2}\right), 
SW=\left(-\frac{1}{\sqrt2},-\frac{1}{\sqrt2}\right),
SE=\left(\frac{1}{\sqrt2},-\frac{1}{\sqrt2}\right).
$
Just as in virtual link diagrams, each double point of a virtual tangle diagram is either a classical or virtual crossing. 

\begin{figure}
\centering
\includegraphics[width=0.4\linewidth]{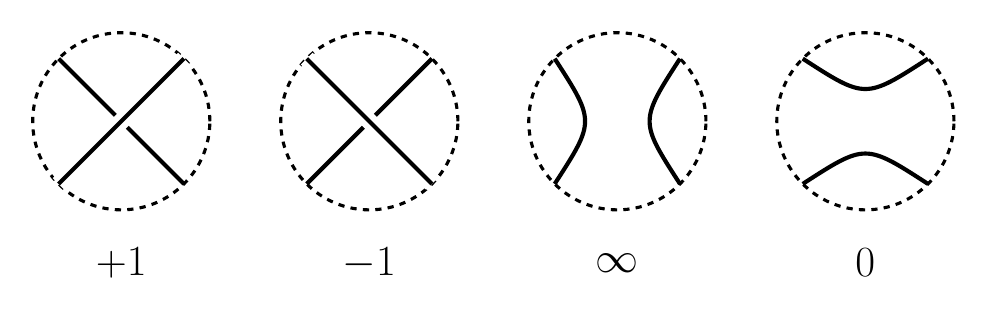}  
\caption{}
\label{fig6}
\end{figure}

Two virtual tangle diagrams are said to be \textit{equivalent} if they are related by a finite sequence of planar isotopies fixing $\partial D^2$ and generalized Reidemeister moves in the interior of $D^2$. A \textit{virtual tangle} is an equivalence class of virtual tangle diagrams. A \textit{classical tangle} is a virtual tangle admitting a diagram with no virtual crossings. Figure \ref{fig6} depicts the elementary rational tangles denoted by $+1$, $-1$, $\infty$, and $0$. 

A virtual tangle diagram $T$ is said to be \textit{disconnected} if either there exists a simple closed curve embedded in $D^2$ disjoint from $T$ but encircles a part of it, or there exists a simple arc properly embedded in $D^2$ which does not meet $T$ and splits $D^2$ into two disks each one containing a part of $T$. A virtual tangle diagram is \textit{connected} if it is not disconnected. A virtual tangle diagram is \textit{reduced} if its number of crossings cannot be decreased by any tangle equivalence. In the remaining of the paper, all virtual tangle diagrams are assumed to be connected and reduced. A virtual tangle diagram is \textit{alternating} if the crossings alternate between over- and under-crossing as one travels along each component. A virtual tangle is said to be \textit{alternating} if it admits an alternating diagram. 

\begin{figure}
\centering
\includegraphics[width=0.3\linewidth]{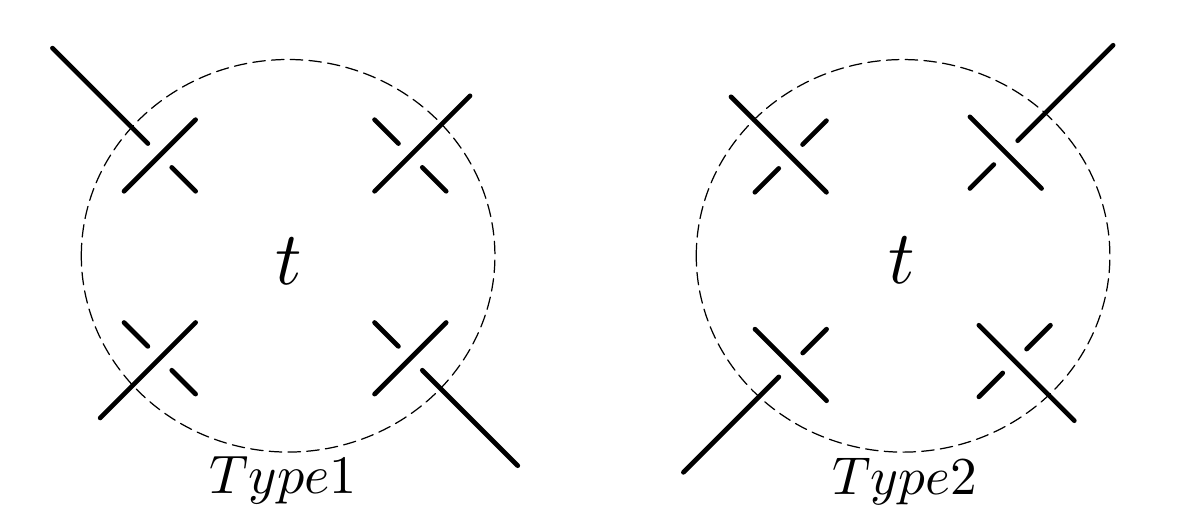}   
\caption{}
\label{fig3}
\end{figure}

For an alternating virtual tangle diagram $T$, it is easy to see that the long arcs of $T$ coming from diametrically opposite endpoints both pass over or below at the first encountered crossing. That remark enables us to distinguish two types of alternating virtual tangle diagrams which we call \textit{type 1} tangles and \textit{type 2} tangles as shown in the Figure \ref{fig3}. 

\begin{figure}
\centering
\includegraphics[width=0.4\linewidth]{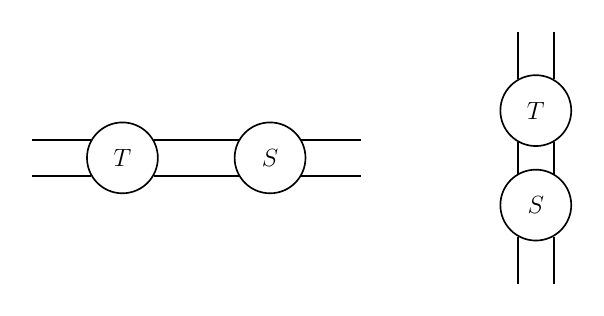}   
\caption{}
\label{tikz8}
\end{figure}

Given two virtual tangle diagrams $T$ and $S$, one can construct two new virtual tangle diagrams $T + S$ and $T * S$ by horizontal and vertical concatenation, respectively, as depicted in Figure \ref{tikz8}. 

\begin{figure}
\centering
\includegraphics[width=0.25\linewidth]{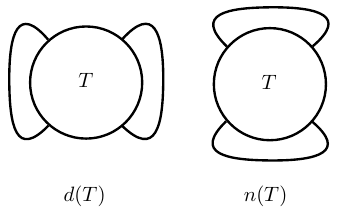}   
\caption{}
\label{fig2}
\end{figure}

A virtual tangle diagram $T$ provides two virtual link diagrams: the \textit{numerator} of $T$, denoted by $n(T)$, which is obtained by joining with simple arcs the two upper endpoints $(NW,NE)$ and the two lower endpoints $(SW,SE)$ of $T$, and the \textit{denominator} of $T$, denoted by $d(T)$, which is obtained by joining with simple arcs each pair of the corresponding top and bottom endpoints $(NW,SW)$ and $(NE,SE)$ of $T$ (see Figure \ref{fig2}). 

\begin{figure}
\centering
\includegraphics[width=0.5\linewidth]{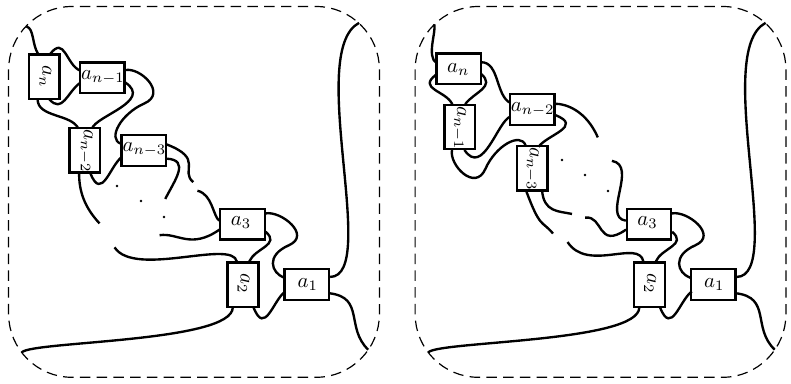}  
\caption{The rational tangle diagram associated to the continued fraction $\left[ a_1,\cdots,a_{n} \right]$ according to whether $n$ is even (left) or odd (right).}
\label{fig7}
\end{figure}

A unique connected classical tangle diagram, as shown in Figure \ref{fig7}, can be associated to each continued fraction $\left[ a_1, a_{2}, \cdots,a_{n} \right]:=\displaystyle a_1 + \cfrac{1}{ a_2  + \cfrac{1}{  \ddots +\cfrac{1}{ a_n  }}},\,\ a_i\in\mathbb{Z}$, where the tangle labeled $a_i$ is the sum of $\abs{a_i}$ copies of the elementary tangle $\pm1$ depending on the sign of $a_i$. It is clear that if all the integers $a_i$ are of the same sign, then the associated tangle diagram is alternating. The equivalence class of such a tangle is called a \textit{rational tangle}. 

In \cite{conway}, Conway showed that for any rational tangle $t$, there exists a unique continued fraction $\left[ a_1  , \cdots ,  a_n \right]$ where all the integers $a_i$ are of the same sign, called the \textit{slope} of $t$, such that the associated alternating tangle diagram depicted in Figure \ref{fig7} represents $t$. Furthermore, Conway showed that two rational tangles are equivalent if and only if they have the same slope. Thus, any rational tangle $t$ can be represented by a unique rational number $\left[ a_1  , \cdots ,  a_n \right] = \frac{p}{q}$ where all the $a_i$ are of the same sign. We note that absolute value of the slope of a rational tangle $t$ is the same as the ratio $\frac{\det(n(t))}{\det(d(t))}$ (see \cite{conway}).

A number  $\varepsilon(c) \in \left\lbrace -1,1 \right\rbrace$ can be assigned to any crossing $c$ of any link diagram according to whether this crossing represents the elementary rational tangle diagram $-1$ or $+1$ as in Figure \ref{fig6}, respectively. The twisting technique introduced for the first time by Champanerkar and Kofman in \cite{kofman2009} replaces such a crossing by any rational tangle $t$ that extends it. More precisely, we say that a rational tangle $t=\left[ a_1, ... , a_n \right]$ \textit{extends} the crossing $c$ if $t$ contains $c$ and for all $i$, $a_i . \varepsilon(c) \geq 1$. One can always replace a crossing $c$ in some virtual link diagram $D$ with a rational tangle $t$ which extends it to get a new link diagram which we will denote by $D^{c \leftarrow t}$. This diagrammatic operation is called a \textit{rational extension} of $c$ with $t$. It is clear that if $D$ is checkerboard colorable then $D^{c \leftarrow t}$ remains so also. 

\section{Determinant of checkerboard colorable links}
\label{sec3}
We start this section by recalling two interpretations of the determinant of checkerboard colorable links introduced by Boden and Karimi in \cite{boden2022} and \cite{boden2023}, respectively. Then, we give a new combinatorial formulation of the determinant.
\subsection{Determinant via coloring and mock Seifert matrices}
\label{section2.2}
Let $L$ be a checkerboard colorable link and let $D$ be a diagram of $L$ with crossings $c_1, \hdots, c_n$ and long arcs $a_1, \hdots, a_m$. Define the \textit{coloring matrix} of $D$, denoted by $B(D)$, as the following $n \times m$ matrix: 
$$b_{ij}(D) = \begin{cases} 
2, &\text{ if } a_j \text{ is the over-crossing arc at the crossing } c_i, \\
-1, &\text{ if } a_j \text{ is an under-crossing arc at the crossing } c_i, \\
0, &\text{ otherwise}.
\end{cases}$$

Boden and Karimi \cite{boden2022} showed that the $(n-1) \times (n-1)$ minors of $B(D)$ are all equal up to sign. The absolute value of such $(n-1) \times (n-1)$ minor is called the \textit{determinant} of $D$, denoted by $\det(D)$. Boden and Karimi also showed that is an invariant of unoriented checkerboard colorable links. Then we can define $\det(L)$ for a checkerboard colorable link $L$ as being the determinant computed using any diagram of $L$. One notable fact about the determinant is that it is invariant under mirror image, i.e., the reversal of the height information at every crossing.

The other interpretation of the determinant of a checkerboard colorable link $L$, given by Boden and Karimi in \cite{boden2023}, is expressed in terms of a matrix which is called the \textit{mock Seifert matrix} and is denoted by $S_L$. To construct $S_L$, one must choose an orientation for $L$ and a \textit{spanning surface}, namely, a compact and connected but not necessarily orientable surface that is bounded by $L$ and embedded in the thickened Carter surface. The determinant of $S_L$ depends on both the chosen orientation for $L$ and the $S^{\ast}$-equivalence class of the spanning surface. Recall that the $S^{\ast}$-equivalence is the equivalence relation on surfaces generated by (i) ambient isotopy, (ii) attachment (or removal) of a tube, and (iii) attachment (or removal) of a half-twisted band. Lifting a checkerboard colorable diagram $D$ of $L$ to the link diagram $\mathcal{D}$ in the Carter surface $\Sigma_D$ naturally yields two spanning surfaces for $L$ in $\Sigma_D \times \left[ 0,1 \right]$, namely the \textit{black surface} and the \textit{white surface}, obtained from the black and white regions of the checkerboard coloring of $\Sigma_D \setminus \mathcal{D}$, respectively. By our convention for choosing the checkerboard coloring (Remark \ref{rmrk1}), this terminology is unambiguous. If $g(\Sigma_D) > 0$, then the black and white checkerboard surfaces of a checkerboard colorable link are not $S^{\ast}$-equivalent. However, every spanning surface is $S^{\ast}$-equivalent to one of these checkerboard surfaces \cite[Proposition\,1.6]{boden2022-2}. This implies that there are exactly two $S^{\ast}$-equivalence classes of spanning surfaces for every checkerboard colorable link. It turns out that the determinants of the two mock Seifert surfaces constructed respectively from the black and white surfaces are equal up to sign, yielding the determinant of the link under consideration by taking their respective absolute values. In particular, Boden and Karimi \cite[Theorem\,3.9, Remark\,3.11]{boden2023} show the following

\begin{align}
\abs{\det(S_L)} = \det(L). 
\label{eq01}
\end{align} 

In \cite{hatano2026}, Hatano and Nozaki give a kind of skein relation for the determinant of the mock Seifert matrix that depends on a specific choice of spanning surface and orientation for the considered link. They showed the following formula

\begin{align}
\det(S_D)=\det(S_{D^{c^*}}) - 2\det(S_{D_{0}^{c}}),
\label{eq2}
\end{align}
where the oriented checkerboard colorable diagrams $D$, $D^{c^*}$, and $D_{0}^{c}$ and the corresponding spanning surfaces are exactly as in Figure \ref{sq}.

The following facts about the determinant can be obtained by Proposition 3.2 and the proof of Theorem 4.6 in \cite{boden2022}. 

\begin{align}
\text{If } D \text{ is a split checkerboard colorable link diagram, then } \det(D)= 0. \tag*{\textbf{Fact 1}} 
\label{eq0}
\end{align}

\begin{align}
\text{If } D \text{ is a non-split alternating virtual link diagram, then } \det(D)= \# \mathcal{S}_D, \tag*{\textbf{Fact 2}}
\label{eq1}
\end{align} 
where $\mathcal{S}_D$ is the set of \textit{connected states} of the checkerboard colorable diagram $D$, that is the set of states of $D$ with a single component.
\subsection{Connected state expansion of the determinant}
\label{section3.2}
\begin{figure}
    \centering
\includegraphics[width=0.7\linewidth]{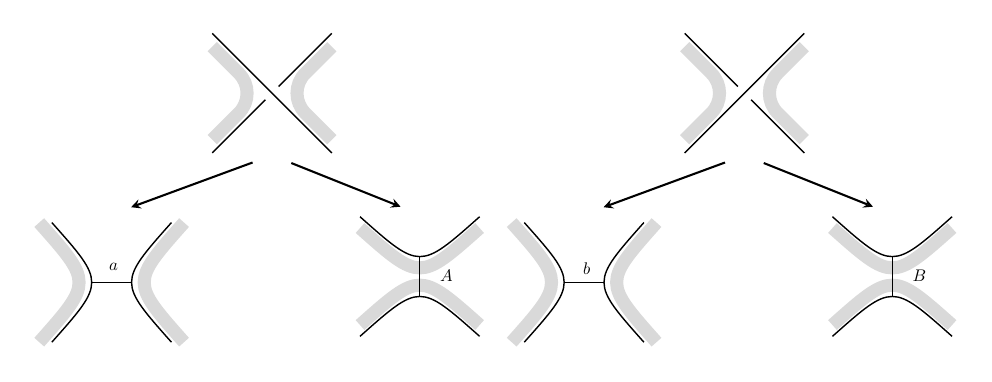}
\caption{At the left (resp. right): crossing with positive (resp. negative) Conway sign and its two smoothings.}
\label{cs}
\end{figure}

\begin{figure}
\centering
\includegraphics[width=0.5\linewidth]{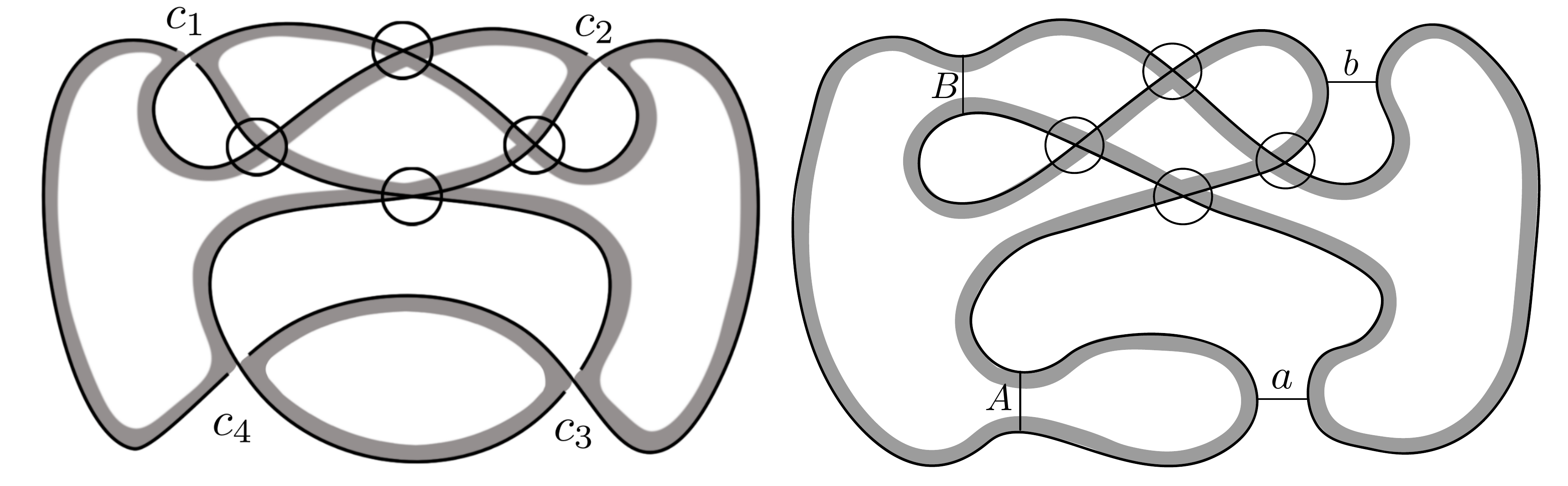}
\caption{}
\label{ccdcs}
\end{figure}

  In each connected state of a checkerboard colorable link diagram, a \textit{marker} is placed in the vicinity of each smoothed crossing, that is a segment connecting the parallel strands resulting from the smoothing. Each marker is labelled by a letter in the alphabet $\left\lbrace A,a,B,b \right\rbrace$ depending on the \textit{Conway sign} of the corresponding crossing as explained in Figure \ref{cs}. In Figure \ref{ccdcs}, we depict a checkerboard colorable link diagram and one of its connected states.

\begin{figure}
\centering
\includegraphics[width=0.6\linewidth]{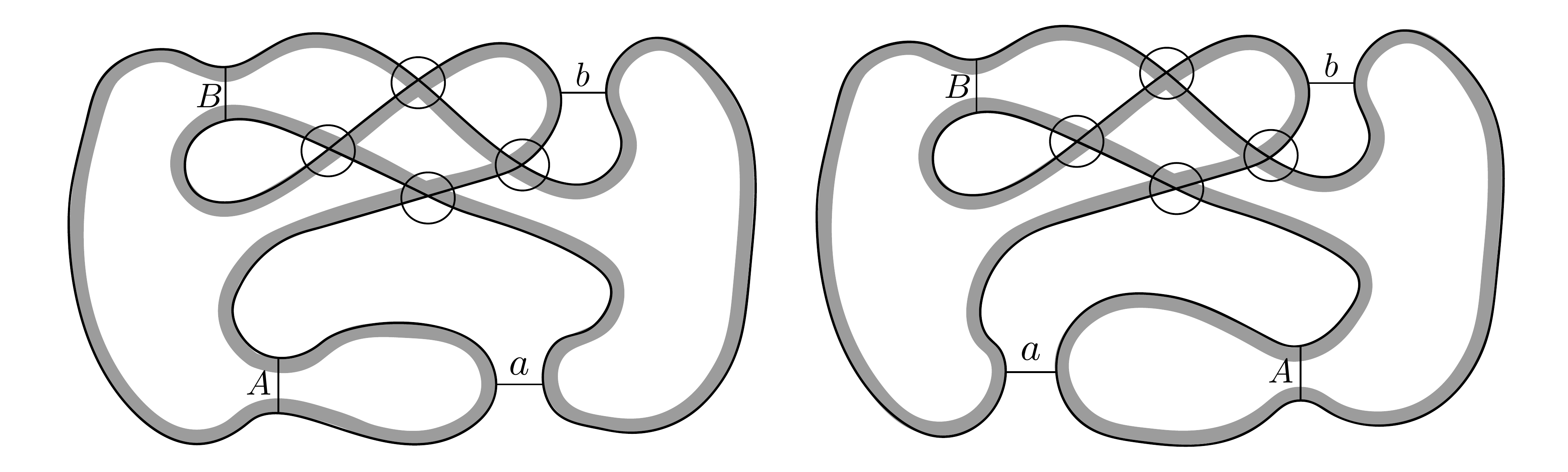}
\caption{}
\label{tr}
\end{figure}

For a connected state $s$ of a checkerboard colorable link diagram $D$, we say that two distinct crossings $c_i$ and $c_j$ \textit{interlace} in $s$ if, when traversing the circle generated by the state $s$, one encounters alternately the smoothing arcs of $c_i$ and $c_j$. For instance, in the connected state on the right of Figure~\ref{ccdcs}, one can easily check that the crossings $c_3$ and $c_4$ interlace, whereas $c_4$ does not interlace with either $c_1$ or $c_2$. By changing the smoothings of an interlacing pair of crossings $c_i$ and $c_j$ in $s$, we obtain a new connected state, denoted $s^{c_ic_j}$. This new state is called the \textit{transposition} of $s$ along the interlacing pair $c_i$ and $c_j$. In Figure~\ref{tr}, we depict the transposition of the connected state on the right of Figure~\ref{ccdcs} along the interlacing crossings $c_3$ and $c_4$. According to \cite[Lemma\,4]{arratia}, connected states form a single orbit under transpositions on interlaced pairs.

If $D$ is an oriented checkerboard colorable diagram with a crossing $c$, we denote by $D^{c}_{\infty}$ the checkerboard colorable link diagram obtained by smoothing $D$ at $c$ that does not inherit orientation from $D$. Unless otherwise specified, $D^{c}_{\infty}$ is considered as an unoriented diagram. Recall that when we consider $D$ as an unoriented diagram, we adopt notations for the diagrams resulting from the smoothing of $c$ as $D^{c}_{\alpha}$ and $D^{c}_{\beta}$ as depicted in Figure \ref{fig.1}. 

\begin{figure}
\centering
\includegraphics[width=0.8\linewidth]{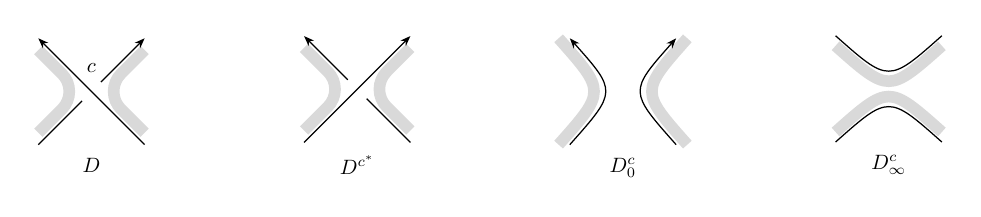}
\caption{}
\label{sq}
\end{figure}

 A crossing $c$ in an oriented checkerboard colorable diagram is \textit{negative} (resp. \textit{positive}) if $c$ is as in the diagram $D$ (resp. $D^{c^*}$) depicted in Figure \ref{sq}. A crossing $c$ of $D$ is said to be \textit{nugatory} if there exists an orientation of $D$ such that its oriented smoothing at $c$ produces a split diagram. If $c$ is a non-nugatory crossing of a checkerboard colorable link diagram $D$, then there exists a one to one correspondence
\[
\mathcal{S}_D \;\longleftrightarrow\; \mathcal{S}_{D^c_{\alpha}} \,\sqcup\, \mathcal{S}_{D^c_{\beta}},
\]
where a state $s \in \mathcal{S}_D$ is mapped to the unique state
$
s' \in \mathcal{S}_{D^c_{\alpha}} \,\sqcup\, \mathcal{S}_{D^c_{\beta}}
$
obtained from $s$ by deleting the marker at the crossing $c$, while keeping all other markers unchanged.
Thus $\# \mathcal{S}_D = \# \mathcal{S}_{D^c_{\alpha}} + \# \mathcal{S}_{D^c_{\beta}}$. Using \ref{eq1}, we get the following fact:

\begin{equation}
\begin{gathered}
\text{If } D \text{ is an alternating virtual link diagram with a non-nugatory crossing } c,
\text{ then}\\[3pt]
\det(D) = \det(D^c_{\alpha}) + \det(D^c_{\beta}).
\end{gathered}
\tag*{\textbf{Fact 3}}
\label{eq5}
\end{equation}

\begin{remark}\label{rmrk1}
For a checkerboard colorable link diagram $D$, there exists a partition of the crossings of $D$ as those with positive Conway sign and those with negative Conway sign as a result of $D$ being alternable and the crossings of an alternating link diagram having the same Conway sign. A checkerboard coloring can be chosen such that the number of crossings with positive Conway sign is maximal. Unless otherwise stated, we shall always choose the checkerboard coloring according to this convention. Therefore, all crossings of an alternating virtual link diagram will have a positive Conway sign. 
\end{remark}

\begin{notation}
For a connected state $s$ of a checkerboard colorable diagram and $l \in \left\lbrace A,a,B,b \right\rbrace$, we denote by $l_s$ the number of markers in $s$ with label $l$.
\end{notation}

\begin{definition}
For a checkerboard colorable link diagram $D$, we define $\delta(D) \in \mathbb{Z}$ as
 \[
\delta(D)=
\begin{cases}
  \sum_{s \in \mathcal{S}_D} (-1)^{a_s},& \text{if}\ D\ \text{is non-split},\\
  0,& \text{if}\ D\ \text{is split}.
\end{cases}
\]
\end{definition}

To prove Theorem \ref{th1}, we show that $\det(D) = \abs{\delta(D)}$ for any checkerboard colorable link diagram $D$. To do so, we organize the argument into several lemmas. We first introduce the following notation that will be used for this section.

\begin{notation} For an oriented link diagram $D$, we denote $D^{i_1 \hdots i_k}_{0_{j_1 \hdots j_m}}$ (resp. $D^{i_1 \hdots i_k}_{\infty_{j_1 \hdots j_m}}$) the link diagram obtained from $D$ by changing the crossings $c_{i_1} \hdots c_{i_k}$ and perform the oriented (resp. unoriented) smoothing at the crossings $c_{j_1} \hdots c_{j_m}$. We also introduce the notations $D^{i_1 \hdots i_k}:=D^{i_1 \hdots i_k}_{\emptyset}$ and $D_{\omega_{j_1 \hdots j_m}}:=D^{\emptyset}_{\omega_{j_1 \hdots j_m}}$, where $\omega \in \left\lbrace 0,\infty \right\rbrace$. 
\end{notation}

\begin{lemma}
Let $D$ be an oriented alternating virtual link diagram with a fixed non-nugatory negative crossing $c$ with positive Conway sign. We let $d := D^{c \leftarrow \frac{1}{2}}$ and we label the crossings of the extension tangle $\frac{1}{2}$ so that $c_1$  and $c_2$  are the top  and bottom crossings, respectively. Then, we have the following 
\begin{itemize}
\item[(i)] $\det(S_{d^1}) = - \det(S_{D^{c}_0}),$
\item[(ii)] $\det ( S_D ) \cdot \det ( S_{D_{0}^{c}} ) < 0.$
\label{lem1}
\end{itemize}
\end{lemma}

\begin{proof}

\item[(i)] It is clear that $d_{0_1} = D$, $d_{0_{12}}=D^{c}_0$, $d^{2}_{0_1}=D^{c^*}$. Using these identifications and the skein relation in Equation (\ref{eq2}), we obtain

\begin{align}
\det(S_d) = \det(S_{d^1})-2\det(S_D)
\label{eq3}
\end{align}

\begin{align}
\det(S_D) = \det(S_{D^{c^*}})-2\det(S_{D^c_0})
\label{eq4}
\end{align}

It is clear that $d_{0_1,\infty_2}\sim d_{\infty_1} \sim D^{c}_{\infty}$. Since these diagrams are all alternating \ref{eq1} and \ref{eq5} imply that

\begin{align}
\det(d)&=\det(d_{0_1})+\det(d_{\infty_1}) \nonumber \\
&= \det(D) + \det(D^{c}_{\infty}) \nonumber \\
&=(\det(D^c_0)+\det(D^c_{\infty}))+\det(D^{c}_{\infty}) \nonumber \\
&=\det(D^c_0)+2\det(D^c_{\infty}).
\label{eq6}
\end{align}

As a result of $d^1 \sim D^{c}_0$, Equation (\ref{eq01}) implies that $\abs{\det(S_{d^1})} = \abs{\det(S_{D^{c}_0})}$. Suppose on the contrary that $\det(S_{d^1}) = \det(S_{D^{c}_0})$. Using \ref{eq1}, Equation (\ref{eq3}), and the fact that $d$, $D$, and  $D^{c}_0$ are alternating diagrams, we obtain

\begin{align}
-\det(S_d) & = 2 \det(S_D) - \det(S_{D^{c}_0}) \nonumber \\
\Rightarrow \abs{\det(S_d)} & = \abs{2 \det(S_D) - \det(S_{D^{c}_0})} \nonumber \\
\xLeftrightarrow{\mathit{(\ref{eq01})}} \det(d) & =  \abs{2 \det(S_D) - \det(S_{D^{c}_0})} \nonumber \\
\xLeftrightarrow{\mathit{(\ref{eq6})}} \det(D^c_0)+2\det(D^c_{\infty}) & = \abs{2 \det(S_D) - \det(S_{D^{c}_0})} \nonumber \\
\xLeftrightarrow{\mathit{\ref{eq5}}} 2\det(D)-\det(D^c_0) & = \abs{2 \det(S_D) - \det(S_{D^{c}_0})} \nonumber \\
\iff \abs{2 \det(S_D)} - \abs{\det(S_{D^{c}_0})} & = \abs{2 \det(S_D) - \det(S_{D^{c}_0})}. \nonumber
\end{align}
It follows that $\det(S_D) \cdot \det(S_{D^{c}_0}) > 0 \nonumber$, which yields the following
\begin{align}
\abs{\det(S_D) + 2\det(S_{D^{c}_0})} & = \abs{\det(S_D)} + 2\abs{\det(S_{D^{c}_0})} \nonumber  \\
\xLeftrightarrow{\mathit{(\ref{eq4})}} \abs{\det(S_{D^{c^*}})} & = \abs{\det(S_D)} + 2\abs{\det(S_{D^{c}_0})} \nonumber  \\
\iff \det(D^{c^*}) & = \det(D)+2\det(D^c_0) \nonumber  \\ 
\xLeftrightarrow{\mathit{\ref{eq5}}} \det(D^{c^*}) & = 3\det(D^c_0) + \det(D^c_{\infty}). 
\label{eq7}
\end{align}

If $D$ is a classical link diagram, then Equation (\ref{eq7}) and Proposition 3.1 in \cite{abchir2020} yield $3\det(D^c_0) + \det(D^c_{\infty})=\abs{\det(D^c_0) - \det(D^c_{\infty})}$. This last equality is known not to hold for a non-nugatory crossing in the classical setting which implies that the invariant under considerartion is not a generalization of the determinant of classical links, which is not the case.

\item[(ii)] By Equation (\ref{eq3}) and (i), we get the following 

\begin{align*}
-2\det(S_D) & = \det(S_{D^c_0})+\det(S_d) \\
\xRightarrow{\mathit{(\ref{eq01})}} 2\det(D) & = \abs{\det(S_{D^c_0})+\det(S_d)} \\
\xLeftrightarrow{\mathit{\ref{eq5}}} 2(\det(D^c_0) + \det(D^c_{\infty})) & = \abs{\det(S_{D^c_0})+\det(S_d)} \\ 
\xLeftrightarrow{\mathit{(\ref{eq6})}} \det(D^c_0) + \det(d) & = \abs{\det(S_{D^c_0})+\det(S_d)} \\
\xRightarrow{\mathit{(\ref{eq01})}} \abs{\det(S_{D^c_0})}+\abs{\det(S_d)} & = \abs{\det(S_{D^c_0})+\det(S_d)} \\
\iff \det(S_d) \cdot \det(D^c_0) & > 0 \\
\xRightarrow{\mathit{(\ref{eq3}),(i)}} \det ( S_D ) \cdot \det ( S_{D_{0}^{c}} ) & < 0.
\end{align*}

\end{proof}

\begin{lemma}
\begin{itemize}
\item[(i)] If $D$ is an oriented virtual link diagram with a negative crossing $c$ of positive Conway sign, then
\begin{equation} \delta(D) = \delta(D^{c^*}) - 2\delta(D^{c}_0). 
\label{eq11}
\end{equation}
\item[(ii)] If $D$ is a non-split alternating virtual link diagram, then there exists an integer $a_D \in \mathbb{Z}_2$, depending only on $D$, such that 
\begin{equation} \delta(D) = (-1)^{a_D} \# \mathcal{S}_D. 
\label{eq12}
\end{equation}
\end{itemize}
\label{lemma3.4}
\end{lemma}

\begin{proof}
\item[(i)] 
We denote $\mathcal{S}_{c,l}(D)$ to be the set of connected states of $D$ where the label of the marker corresponding to the crossing $c$ is $l\in \left\lbrace A,a,B,b \right\rbrace$. As a consequence of $c$ being negative and of positive Conway sign, we have a natural one-to-one correspondence between the sets $\mathcal{S}_{c,a}(D)$, $\mathcal{S}_{c,b}(D^{c^*})$, and $\mathcal{S}_{D^c_0}$. Analogously, we have a natural one-to-one correspondence between the sets $\mathcal{S}_{c,A}(D)$, $\mathcal{S}_{c,B}(D^{c^*})$, and $\mathcal{S}_{D^c_{\infty}}$. For $s_1 \in \mathcal{S}_{c,a}(D)$ and $s_2 \in \mathcal{S}_{c,A}(D)$, we denote $s'_1 \in \mathcal{S}_{c,b}(D^{c^*})$ and $s'_2 \in \mathcal{S}_{c,B}(D^{c^*})$ the corresponding states of $s_1$ and $s_2$, respectively. We also denote $\sigma_1$ and $\sigma_2$  the states of $D^c_0$ and $D^c_{\infty}$ corresponding to $s_1$ and $s_2$, respectively. It is clear that $a_{s_1}=a_{s'_1}+1=a_{\sigma_1}+1$, and  $a_{s_2}=a_{s'_2}=a_{\sigma_2}$. This allows us to express $\delta(D)$ and $\delta(D^{c^*})$ in terms of $\delta(D^c_0)$ and $\delta(D^c_{\infty})$ as follows

\begin{align}
\delta(D) &=\sum_{s \in \mathcal{S}_D}(-1)^{a_s} \nonumber  \\
&= \sum_{s \in \mathcal{S}_{c,a}(D)}(-1)^{a_s} + \sum_{s \in \mathcal{S}_{c,A}(D)}(-1)^{a_s} \nonumber  \\
&= \sum_{\sigma \in \mathcal{S}_{D^c_0}}(-1)^{a_{\sigma}+1} + \sum_{\sigma \in \mathcal{S}_{D^c_{\infty}}}(-1)^{a_{\sigma}} \nonumber  \\
\delta(D) &= -\delta(D^c_0) + \delta(D^c_{\infty}).
\label{sr1}
\end{align}

\begin{align}
\delta(D^{c^*}) &=\sum_{s' \in \mathcal{S}_{D^{c^*}}}(-1)^{a_{s'}} \nonumber  \\
&= \sum_{s' \in \mathcal{S}_{c,b}(D^{c^*})}(-1)^{a_{s'}} + \sum_{s' \in \mathcal{S}_{c,B}(D^{c^*})}(-1)^{a_{s'}} \nonumber  \\
&= \sum_{\sigma \in \mathcal{S}_{D^c_0}}(-1)^{a_{\sigma}} + \sum_{\sigma \in \mathcal{S}_{D^c_{\infty}}}(-1)^{a_{\sigma}} \nonumber  \\
\delta(D^{c^*}) &= \delta(D^c_0) + \delta(D^c_{\infty}).
\label{sr2}
\end{align}

The result then follows by subtracting Equations (\ref{sr2}) and (\ref{sr1}).

\item[(ii)] Let $D$ be a non-split alternating virtual link diagram. Since all crossings of $D$ are of positive Conway sign, the labels of the markers in all the states are either $A$ or $a$. From Table \ref{tab1}, we see that $a_s$ is invariant modulo $2$ under any transposition of $s$. Therefore, since connected states form a single orbit under transpositions on interlaced pairs (\cite[Lemma\,4]{arratia}), we define unambiguously the number $a_D \in \mathbb{Z}_2$ and obtain $\delta(D) = \sum_{s \in \mathcal{S}_D} (-1)^{a_s} = \sum_{s \in \mathcal{S}_D} (-1)^{a_D} = (-1)^{a_D} \# \mathcal{S}_D$. 
\end{proof}

\begin{remark}
In contrast to the classical case, not all virtual knots (or links) can be unknotted (or unlinked) by crossing changes. Indeed, crossing changes preserve the homotopy class, so homotopically nontrivial virtual knots (or links) cannot be trivialized by these moves \cite[Section\,1.2]{boden2019}. Therefore, we can not only compare the recursive formulas of $\det(S_D)$ and $\delta(D)$ in Equations (\ref{eq2}) and (\ref{eq11}) to conclude the result in Theorem \ref{th1} as the stop condition is not always attainable. 
\end{remark}

\begin{figure}
\centering
\includegraphics[width=0.4\linewidth]{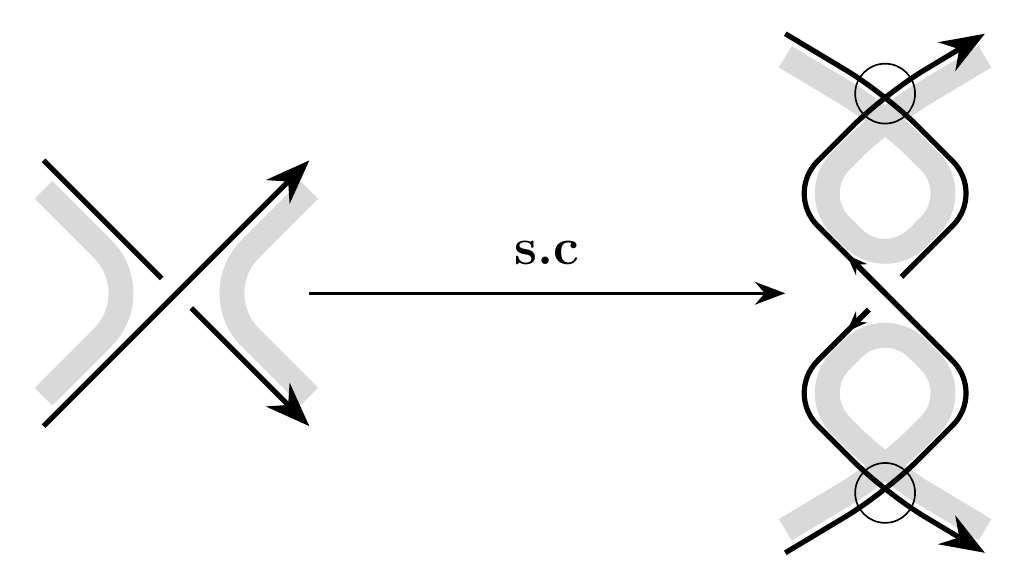}
\caption{}
\label{sc}
\end{figure}

\begin{proof}[Proof of Theorem] \ref{th1}

Let $C_-(D) = \lbrace c_1, ... , c_k \rbrace$ be the set of crossings with negative Conway sign in the oriented checkerboard colorable link diagram $D$. We perform the \textit{sign change} ($\textbf{s.c}$) operation at each negative crossing in $C_-(D)$ as depicted in Figure \ref{sc}. We denote $\D$ the obtained checkerboard colorable link diagram. The coloring matrices of $D$ and $\D$ are clearly equal, and thus, $\det(D)=\det(\D)$.  Since the number $\delta(D)$ is invariant under the sign change operation when applied at a negative crossing of negative Conway sign, $\delta(D)=\delta(\D)$. To show that $\det(D) = \abs{\delta(D)}$, we only need to verify that $\abs{\delta(\D)} = \abs{\det(S_{\D})}$.

To compute $\det(S_{\D})$ using the recursive relation in Equation (\ref{eq2}), we construct a binary tree which is called an \textit{alternation resolution tree}, as follows: We start with the root diagram $\D$ and perform a crossing change and an oriented smoothing at the crossing $c_1$. Then we connect $\D$ to $\D_{0_1}$ and $\D^1$ by drawing a line segment pointing to the bottom left and right to which we assign the label $2$ and $1$, respectively. We then repeat this procedure for $\D_{0_1}$ and $\D^1$ using the crossing $c_2$, and so on, until the alternating diagrams $\D_{0_{1\cdots k}}$ and $\D^{1\cdots k}$ appear as leaves of the resulting tree (see Figure \ref{rt}). If we denote $d_1,\ldots,d_{2^k}$ the leaves of this tree, ordered from right to left, as in Figure \ref{rt}, then we can write 

\begin{align}
\det(S_{\D}) = \sum_{i} \eta_i \det(S_{d_i}), 
\label{eq9}
\end{align}
where the integer $\eta_i$ is obtained by multiplying the coefficients on the line segments along the uniquely determined branch path that begins with the root $\D$ and terminates at the leave $d_i$. We observe that $d_2, \hdots, d_{2^k}$ are all alternating diagrams obtained from the alternating diagram $d_1$ by oriented smoothings. Some of the diagrams $d_2, \hdots, d_{2^k}$ may be split after some index, and hence, by \ref{eq0}, will not contribute to the sum in Equation (\ref{eq9}). Thus, by \ref{eq1}, we have the following

\begin{align}
\det(S_{\D}) = \sum_{i \text{ such that } d_i \text{ is non-split }} (-1)^{\epsilon_i} \eta_i \# \mathcal{S}_{d_i}, 
\label{eq10}
\end{align}
where the integer $\epsilon_i \in \mathbb{Z}_2$ is defined as $\epsilon_i = \epsilon_1$ if $d_i$ is obtained from $d_1$ by an even number of smoothings, or $\epsilon_i = \epsilon_1 + 1$ if $d_i$ is obtained from $d_1$ by an odd number of smoothings as a direct consequence of Lemma \ref{lem1}. 

Using the same alternation resolution tree, the result in Equation (\ref{eq12}), and the recursive formula in Equation (\ref{eq11}), we can analogously write

\begin{align}
\delta(\D) = \sum_{i \text{ such that } d_i \text{ is non-split }} (-1)^{a_{d_i}} \eta_i \# \mathcal{S}_{d_i}.
\label{eq13}
\end{align}

Since the crossings $c_1, \hdots, c_n$ are negative and of positive Conway sign in the diagram $d_1$, we have $a_{d_i} = a_{d_1}$ if $d_i$ is obtained from $d_1$ by an even number of smoothings, and $a_{d_i} = a_{d_1} + 1$ if $d_i$ is obtained from $d_1$ by an odd number of smoothings. Thus, Equations (\ref{eq10}) and (\ref{eq13}) imply that
$$
 \delta(\D) =
\begin{cases}
  \det(S_{\D}), & \text{if}\ \  a_{d_1} = \epsilon_1,\\ - \det(S_{\D}),& \text{if}\ \  a_{d_1} = \epsilon_1 + 1. 
\end{cases}
$$
Then we get $\abs{\delta(\D)} = \abs{\det(S_{\D})}$. This completes the proof.
\begin{figure}
\centering
\includegraphics[width=1\linewidth]{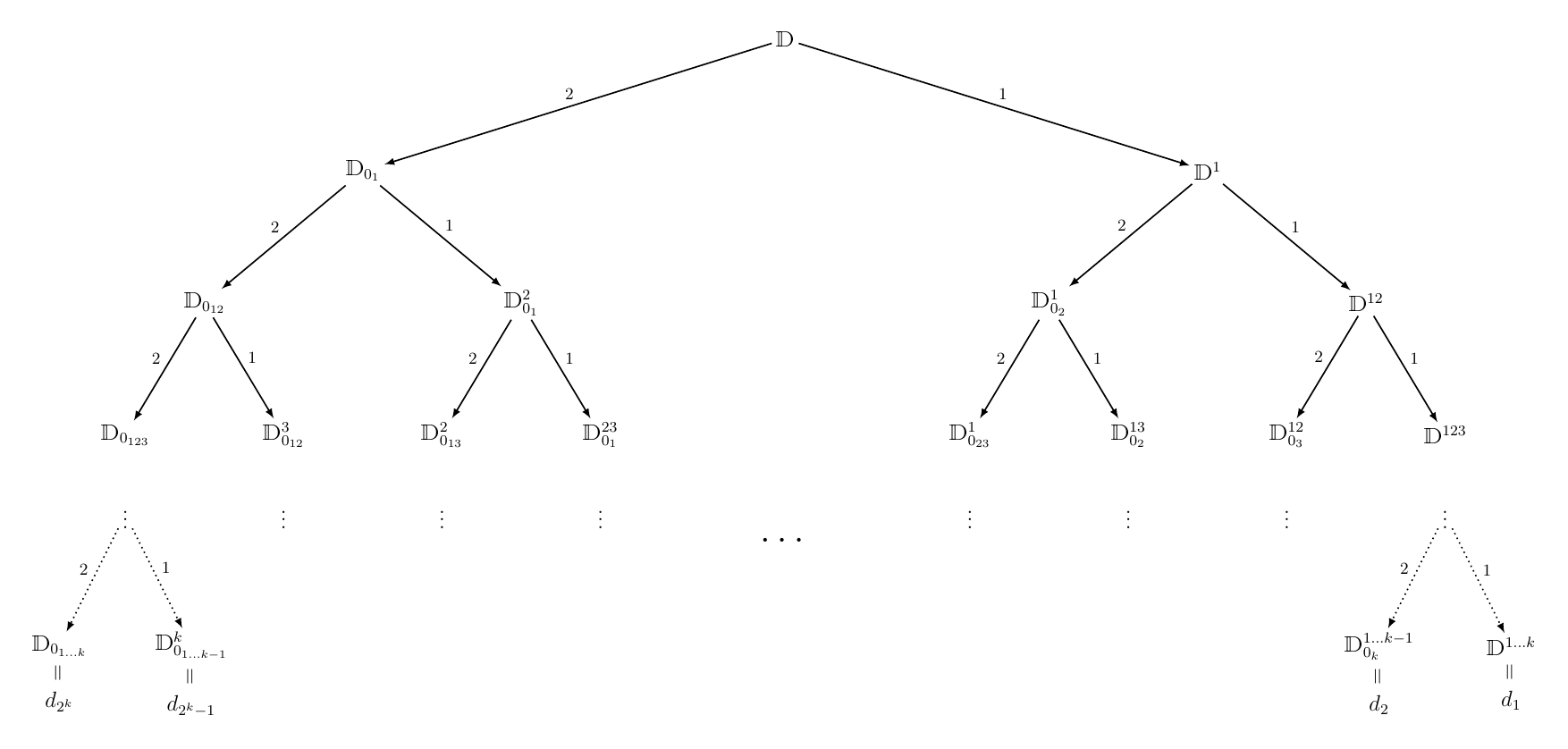}
\caption{The alternation resolution tree of the diagram $\D$.}
\label{rt}
\end{figure}

\end{proof}

\begin{example} 
In Figure \ref{diagramplusstates}, we depict a checkerboard colorable knot diagram that represents the knot $4.99$ in Green's table of virtual knots \cite{green}, together with its five connected states. By Theorem \ref{th1}, we get $\det(4.99) = 3$. 

\begin{figure}
\centering
\includegraphics[width=1\linewidth]{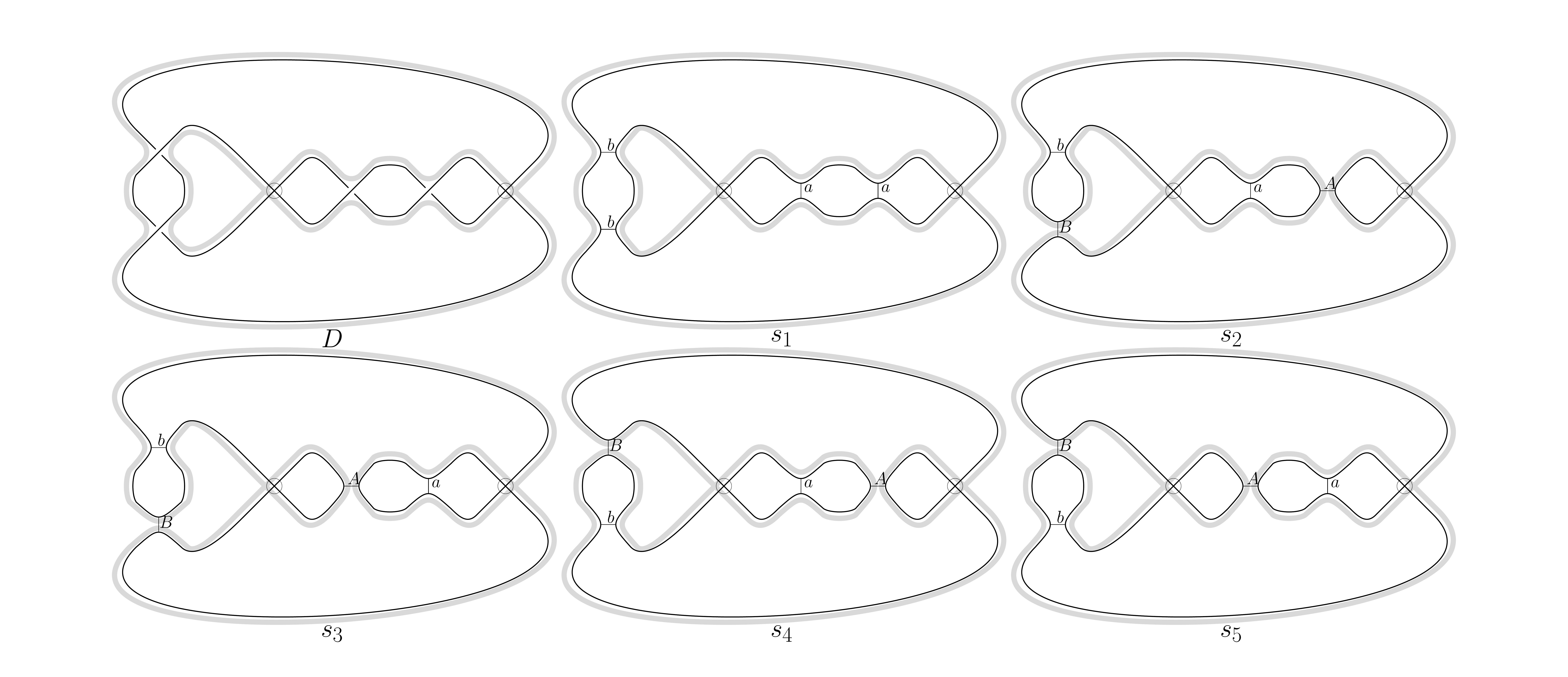}   
\caption{}
\label{diagramplusstates}
\end{figure}
\label{expl3.5}
\end{example}

Now we apply the connected state expansion of the determinant to extend some known determinant formulas in the classical setting to the virtual one. The following result is a key tool in establishing these formulas.

\begin{proposition}
Let $c$ be a crossing of a checkerboard colorable link diagram such that $\det(D) = \det(D^c_{\alpha}) + \det(D^c_{\beta})$ and let $\dfrac{p}{q}$ be a rational number. Then, we have the following
\[
  \det \left( D^{c \leftarrow \frac{p}{q}} \right) =
  \begin{cases}
     p \det(D^c_{\alpha}) + q \det(D^c_{\beta}), & \text{if }\ \varepsilon(c) = +1, \\
     p \det(D^c_{\beta}) + q  \det(D^c_{\alpha}), & \text{if }\ \varepsilon(c) = -1.
  \end{cases}
\]
\label{prop3.11}
\end{proposition}

\begin{proof}
We prove the first case for which $\varepsilon(c) = +1$, the other case can be established by a similar argument. We may assume that the crossing $c$ has a negative Conway sign by changing the checkerboard coloring if needed. We prove the result by induction on the length of the continued fraction of the given rational tangle. If the length is one, then we have an integer tangle. We prove this case by induction on the number of crossings of this integer tangle. Both such inductions depend on proving the result for the diagrams $d:=D^{c \leftarrow \frac{1}{2}}$ and $d' :=D^{c \leftarrow 2}$ where the crossing $c$ has been replaced by a vertical and a horizontal integer tangle of two crossings that extends $c$, respectively. Now we consider the first diagram $d$. If we smooth $d$ at anyone of these two crossings, which will be denoted also by $c$, we obtain the links $d^c_{\alpha}$ and $d^c_{\beta}$ that represent the diagrams $D$ and $D^c_{\beta}$, respectively.  By using Equation (\ref{sr2}) and since $\det(D) = \det(D^c_{\alpha}) + \det(D^c_{\beta})$ yields the inequality $\delta(D^c_{\alpha}) \cdot \delta(D^c_{\beta}) > 0$, we get
  $$\det(d) = |\delta(D^{c \leftarrow \frac{1}{2}})|   = |\delta(D) +\delta(D^c_{\beta}) |  = |  \delta(D^c_{\alpha}) + 2 \delta(D^c_{\beta})|  =  | \delta(D^c_{\alpha}) | + 2| \delta(D^c_{\beta})|,$$ hence
$$\det(d) = \det(d^c_{\alpha})  + \det(d^c_{\beta}) = \det(D) + \det(D^c_{\beta}) = \det(D^c_{\alpha}) + 2\det(D^c_{\beta}).$$

Thus we conclude that the result holds for $D^{c \leftarrow \frac{1}{2}}$. For $d'$, the result can be obtained by a similar argument. Now we assume that the result holds for $D^{c \leftarrow \frac{1}{n-1}}$. We incorporate the previous argument for $D^{c \leftarrow \frac{1}{2}}$ to show that the result holds for $D^{c \leftarrow \frac{1}{n}}$. Similarly, we also obtain that the result holds for $D^{c \leftarrow n}$. Suppose now that the result holds for any rational tangle extending $c$ whose fraction has length $n$. Note that $[a_{1},a_{2}, \ldots, a_{n}+1]=[a_{1},a_{2}, \ldots, a_{n},1]$ and that $[a_{1},a_{2}, \ldots, a_{n+1}]$ is obtained from $[a_{1},a_{2}, \ldots, a_{n},1]$ by extending the last crossing with the integer tangle $a_{n+1}$ which can be vertical or horizontal depending the the parity of $n$ (see Figure \ref{rt}). At the end, by using the first induction on the number of crossings in the integer tangle $a_{n+1}$, we obtain that the result holds for any rational tangle that extends the crossing $c$.
\end{proof}

We now introduce the class of almost alternating virtual links, defined as a natural extension of the classical notion. As we shall see, the determinant formula of classical almost alternating links extends naturally to the virtual setting.

\begin{definition}
An \textit{almost alternating} virtual link diagram $D$ is a non-alternating virtual link diagram that becomes alternating by the change of only one crossing. Such crossing is called a \textit{dealternator} of $D$. A virtual link is \textit{almost alternating} if it is not alternating and admits an almost alternating link diagram.
\end{definition}

Any almost alternating virtual link diagram is clearly alternable and hence checkerboard colorable. In \cite[Proposition\,2.4]{tsukamoto}, Tsukamoto shows that if $D$ is a connected classical almost alternating diagram with more than one dealternator, then $D$ is obtained from the standard Hopf link diagram with two crossings by changing one of them. Tsukamoto's proof depends only on the adjacency of crossings, i.e., how the crossings are "connected" by the arcs in the given diagram. On the other hand, the presence of virtual crossings does not affect that argument. Consequently, the same argument can be applied verbatim to the virtual setting, yielding the following proposition.

\begin{proposition}
Let $D$ be a connected virtual almost alternating link diagram. If $D$ has more than one dealternator, then $D$ is either obtained from the standard Hopf link diagram with two crossings by changing one of them or the virtual link diagram 2-2-1 in Bartholomew's table of virtual links \cite{bartholomew}.
\label{prop1}
\end{proposition}

\begin{remark} From now on, we only consider almost alternating link diagrams having only one dealternator by exluding the two cases mentioned in Proposition \ref{prop1}. We note that both of these links have zero determinant.
\end{remark}

The following result allows us to compute the determinant of an almost alternating link diagram. This formula can be regarded as an extension to the virtual setting of that given in \cite[Proposition\,3.1]{abchir2020}.
\begin{proposition} Let $D$ be an almost alternating link diagram with dealternator $c$. Then
  $$\det(D)=\abs{\det(D^c_{\alpha})-\det(D^c_{\beta})} = \abs{\# \mathcal{S}_{D^c_{\alpha}}-\# \mathcal{S}_{D^c_{\beta}}}.$$
\label{prop2}  
\end{proposition}

\begin{proof} We orient $D$ so that $c$ is a negative crossing. Since $D$ is almost alternating with dealaternator $c$, this crossing has a different Conway sign from other crossings of $D$. We choose the checkerboard coloring so that $c$ has a positive Conway sign. Then, $c$ satisfies the conditions of item (i) in Lemma \ref{lemma3.4}, and then we can write the following identity used in the proof:
$$\delta(D)=-\delta(D_0^c)+\delta(D_\infty^c).$$
We obtain
$$\delta(D)= - \# \mathcal{S}_{D^c_{\alpha}}+\# \mathcal{S}_{D^c_{\beta}}=-\det(D_0^c)+\det(D_\infty^c) ,$$
where the last equality is obtained from \ref{eq1} and the fact that $D^c_{\alpha}$ and $D^c_{\beta}$ are both alternating with crossings of negative Conway sign.
\end{proof}

\begin{remark}
Let $D$ be an almost alternating virtual link diagram with dealternator $c$. Whithout loss of generality, we can assume that $c$ is identified with the rational tangle $-1$ and that $D$ is expressed as the numerator closure of the sum of $c$ and an alternating virtual tangle $T$ of type $1$. Thus, every almost alternating virtual link diagram can be written in the form $D = n(-1 + T)$.	
\label{rmrk2}
\end{remark}

We now give the formula for computing the determinant of virtual links obtained as the closure of the sum of an alternating virtual tangle and a rational tangle. The following result can be regarded as an extension of \cite[Proposition\,3.4]{abchir2020} to the virtual setting. 

\begin{proposition}
Let $T$ be an alternating virtual tangle diagram of type $1$ and let $\dfrac{p}{q} > 0$ be a rational number. We have the following 
\begin{enumerate}
\item[(i)] $\det\left(n\left(\dfrac{p}{q}+T\right)\right) = q \det\left(n(T)\right)+ p \det\left(d(T)\right) = q\# \mathcal{S}_{n(T)} + p \# \mathcal{S}_{d(T)}.$
\item[(ii)]$\det\left(n\left(-\dfrac{p}{q}+T\right)\right) = \abs{ q \det\left(n(T)\right)- p \det\left(d(T)\right) }=\abs{ q  \# \mathcal{S}_{n(T)}- p \# \mathcal{S}_{d(T)} }.$
\end{enumerate}
\label{prop3.10} 
\end{proposition}

\begin{proof}
\begin{enumerate}
\item[(i)] Let $c$ be the leftmost crossing of the alternating virtual link diagram $n(1+T)$ denoted $D$. By connectivity of $T$ and \ref{eq1}, the alternating virtual link diagrams $D^c_{\alpha} = d(T)$ and $D^c_{\beta} = n(T)$ have non-zero determinants equal to the numbers of their connected states. This implies that $c$ is not a nugatory crossing. By \ref{eq5}, we have $\det(D) = \det(D^c_{\alpha}) + \det(D^c_{\beta})$. The result now follows from Proposition \ref{prop3.11} and the fact that $\det\left(n\left(\frac{p}{q}+T\right)\right) = \det\left(D^{c \leftarrow \frac{p}{q}}\right) $. 
\item[(ii)] By rotating the diagram in the plane by $90$ degrees then taking the mirror image if needed, we may assume that $q > p$. Since the rational tangles of slopes $-1 + \frac{q-p}{q}$ and $-\frac{p}{q}$ are equivalent $n\left(-1+\left(\frac{q-p}{q}+T\right)\right)$ is an almost alternating diagram equivalent to $n\left(-\frac{p}{q}+T\right)$. By Proposition \ref{prop2} and item (i) of this proposition, we have
 
\begin{align}
\det\left(n\left(-\frac{p}{q}+T\right)\right)& = \det\left(n\left(-1+\left(\frac{q-p}{q}+T\right)\right)\right)\nonumber\\ & = \abs{\det\left(n\left(\frac{q-p}{q}+T\right)\right) - \det\left(d\left(\frac{q-p}{q}+T\right)\right)}.\nonumber
  \end{align}
    The diagram $d\left(\frac{q-p}{q}+T\right)$ is equivalent to the connected sum of $d\left(\frac{q-p}{q}\right)$ and $d(T)$, thus, by \cite[Theorem\,6.4]{boden2022}, we have $\det\left(d\left(\frac{q-p}{q}+T\right)\right) = \det\left(d\left(\frac{q-p}{q}\right)\right)\det\left(d(T)\right) = q \det\left(d(T)\right)$. By item (i), we have $\det\left(n\left(\frac{q-p}{q}+T\right)\right) = q \det\left(n(T)\right) + (q-p) \det\left(d(T)\right)$. Thus, the result follows.
\end{enumerate} 

\end{proof}

\section{The determinant versus the mock-determinant}
\label{sec4}
In this section, we compare the determinant and the mock-determinant of checkerboard colorable links using their respective connected state expansions given in Theorem \ref{th1} and \cite[Proposition\,3.7]{abchir2026}. To do so, we start by giving some definitions.

For a connected state $s$ of a checkerboard colorable link diagram $D$, a crossing $c$ of $D$ is said to be \textit{internal} to $s$ if the marker assigned to $c$ in $s$ is either an $A$-marker or $B$-marker; otherwise, $c$ is \textit{external} to $s$. A pair of crossings of $D$ is said to be \textit{homogeneous} in $s$ if both are either internal or external in $s$; otherwise, the pair is \textit{heterogeneous} in $s$. A connected state $s$ is said \textit{classical} if every pair of interlacing crossings in $s$ is heterogeneous. Since a connected state of a classical link diagram is embedded in the plane, it is easy to check that it is classical using the Jordan curve theorem.

\begin{lemma}
If a checkerboard colorable link diagram has a classical state, then all of its connected states are classical.
\label{lemma4.1}
\end{lemma}

\begin{proof}
  Let $D$ be a checkerboard colorable link diagram and $s$ be a connected state of $D$. As stated in Section \ref{section3.2}, connected states form a single orbit under transpositions on interlaced pairs (\cite[Lemma\,4]{arratia}). So it suffices to show that a classical state cannot be transformed into a non-classical state through a transposition along an interlaced pair of crossings. For every crossing $c$ of $D$, we denote by $C_c(s)$ the set of crossings that interlaces with $c$ in $s$, and by $\overline{C}_c(s)$ the set $C(D) \setminus (C_c(s) \cup \left\lbrace c \right\rbrace)$, where $C(D)$ is the set of all crossings of $D$. Suppose on the contrary that there exists an interlacing pair $uv$ in the classical state $s$ such that the state $s':=s^{uv}$ is non-classical. Let $kl$ be a homogeneous interlaced pair in $s'$. The pairs $kl$ and $uv$ are different since $uv$ is still heterogeneous in $s'$. Suppose that $\lbrace k,l \rbrace \cap \lbrace u,v \rbrace \neq \emptyset$. Assume without loss of generality that $v=l$, and thus $k \neq u$. By the fact that $k \in C_v(s')$ and the second item of \cite[Lemma\,4.11]{abchir2025}, we have that $ku$ is an interlacing pair in $s$. Furthermore, since $kv$ is homogeneous and $uv$ is heterogeneous in $s'$, $ku$ is heterogeneous in $s'$. This implies that $ku$ is a homogeneous interlacing pair in $s$, which contradicts that $s$ is classical. Hence, $\lbrace k,l \rbrace \cap \lbrace u,v \rbrace = \emptyset$. Suppose now that $\left\lbrace k,l \right\rbrace \subset \overline{C}_u(s) \cap \overline{C}_v(s)$, then by the first item of \cite[Lemma\,4.11]{abchir2025}, we get that $kl$ is a homogeneous interlacing pair in $s$, which again contradicts $s$ being classical. Thus, $\left\lbrace k,l \right\rbrace \subset C_u(s) \cup C_v(s)$. By classicality of $s$, we have $C_u(s) \cap C_v(s) = \emptyset$. Therefore, $(k,l) \in \left( C_u(s) \cap \overline{C}_v(s) \right) \times \left( \overline{C}_u(s) \cap C_v(s) \right)$ or $(k,l) \in \left( \overline{C}_u(s) \cap C_v(s) \right) \times \left( C_u(s) \cap \overline{C}_v(s) \right)$. Assume without loss of generality that $(k,l) \in \left( \overline{C}_u(s) \cap C_v(s) \right) \times \left( C_u(s) \cap \overline{C}_v(s) \right)$. By classicality of $s$, if $u$ is internal in $s$, and thus $v$ is external in $s$, we get that $k$ is internal and $l$ is external in $s$, which makes $kl$ a heterogeneous interlacing pair in $s'$, which is absurd. Similarly, we get a contradiction if $u$ is external in $s$.
\end{proof}

To each connected state $s$ of a link diagram $D$, we assign an integer $\nu_s$ defined by $\nu_s := B_s + a_s$. Also, for each pair of positive integers $m$ and $n$, we define the following set: $$\mathcal{S}_D^{n,m} := \left\lbrace s \in \mathcal{S}_D \mid \nu_s = n \pmod 4, a_s = m \pmod 2 \right\rbrace.$$
From Table \ref{tab1}, it is easy to see that $\nu_s$ is invariant modulo $2$ under any transposition of $s$. This uniquely determines an integer $\nu_D \in \mathbb{Z}_2$. Based on this, we define the following integers
\[
\lambda(D):=
\begin{cases}
\# \mathcal{S}_D^{1,0}-\# \mathcal{S}_D^{3,1}, & \text{if } \nu_D=1,\\
\# \mathcal{S}_D^{0,0}-\# \mathcal{S}_D^{2,1}, & \text{if } \nu_D=0.
\end{cases}
\quad,\quad
\gamma(D):=
\begin{cases}
\# \mathcal{S}_D^{1,1}-\# \mathcal{S}_D^{3,0}, & \text{if } \nu_D=1,\\
\# \mathcal{S}_D^{0,1}-\# \mathcal{S}_D^{2,0}, & \text{if } \nu_D=0.
\end{cases}
\]

$$
\mu(D) := \sum_{s \in \mathcal{S}_D} (-1)^{\left\lfloor \frac{\nu_s}{2} \right\rfloor}.\nonumber
$$

\begin{remark}
  According to \cite[Proposition\,3.7]{abchir2026}, we have $\mdet{(D)}= \abs{\sum_{s \in \mathcal{S}_D} i^{-\nu_s}}$. It is easy to see that $\mdet{(D)} = \abs{\mu(D)}.$
\end{remark}

\begin{proposition}
Let $D$ be a checkerboard colorable link diagram, then
\begin{itemize}
\item[(i)] $\det(D) \equiv \mdet{(D)}  \equiv \# D \pmod 2$, where $\# D$ is the number of components of $D$. 
\item[(ii)] If $D$ has a classical state, then $\det(D) = \mdet{(D)}$.
\item[(iii)]  If $D$ is alternating, then $\mdet(D) \leq \det{(D)}$ and the equality holds if and only if $D$ has a classical state.
\end{itemize}
\label{prop4.1}
\end{proposition}

\begin{proof}
\item[(i)] It is clear that $\delta(D) = \lambda(D) - \gamma(D)$ and $\mu(D) = \lambda(D) + \gamma(D)$. On the other hand, we have that $\mdet{(D)} \equiv \# D \pmod 2$ \cite[Proposition\,3.5]{abchir2026}.

\item[(ii)] For a set $S \subset \mathcal{S}_D$, we define $$\tau_{hom} (S)=\left\lbrace s^{ij} \mid s \in S, ij \text{ is a homogeneous interlacing pair in } s \right\rbrace,$$ $$\tau_{het} (S)=\left\lbrace s^{ij} \mid s \in S, ij \text{ is a heterogeneous interlacing pair in } s \right\rbrace.$$ From Table \ref{tab1}, we have that 
$\tau_{hom}(\mathcal{S}_D^{n,m}) \subset \mathcal{S}_D^{n+2,m} \cup \mathcal{S}_D^{n,m+1}$, and $\tau_{het}(\mathcal{S}_D^{n,m}) \subset \mathcal{S}_D^{n,m} \cup \mathcal{S}_D^{n+2,m+1}$. Thus, if $D$ has a classical state, which implies that all connected states of $D$ are classical (Lemma \ref{lemma4.1}), then either $\lambda(D)$ or $\gamma(D)$ is zero. This implies that $\abs{\delta(D)} = \abs{\mu(D)}$. 

\item[(iii)] Assume that $D$ is alternating. Since all crossing of $D$ have positive Conway sign, $\nu_s = a_s$ for every classical state $s$ of $D$. On the other hand, we showed in the proof of Lemma \ref{lemma3.4} that $a_D = a_s \mod 2  $ is constant over all connected states of $D$. Hence, we have $ \abs{\mu(D)} = \abs{\sum_{s \in \mathcal{S}_D} (-1)^{\left\lfloor \frac{a_s}{2} \right\rfloor}} \leq \# \mathcal{S}_D$. Consequently, we get $\mdet(D) \leq \det(D)$. On the other hand, we have $\mathcal{S}_D^{n,m} = \mathcal{S}_D^{n,a_D}$ for every integers $n$ and $m$. Also, we have 
\begin{align*}
 \mdet(D) = \det(D) \iff & \abs{\sum_{s \in \mathcal{S}_D} (-1)^{\left\lfloor \frac{a_s}{2} \right\rfloor}} = \# \mathcal{S}_D \\
\iff & \nu_D := a_s \mod 4 \text{ is constant over all connected states of } D \\
\iff & \bigcup\limits_{n} \left( \tau_{hom}(\mathcal{S}_D^{n,a_D}) \cup \tau_{het}(\mathcal{S}_D^{n,a_D}) \right) \subset \mathcal{S}_D^{\nu_D,a_D} \\
\iff & \tau_{hom}(\mathcal{S}_D^{n,a_D}) = \emptyset \text{, for every integer } n,
\end{align*} 
where the third equivalence comes from the fact that $\mathcal{S}_D =  \bigcup\limits_{n} \left( \tau_{hom}(\mathcal{S}_D^{n,a_D}) \cup \tau_{het}(\mathcal{S}_D^{n,a_D}) \right)$, since the connected states form a single orbit under transpositions on interlaced pairs (\cite[Lemma\,4]{arratia}), and the last equivalence comes from the fact that $\tau_{hom}(\mathcal{S}_D^{n,a_D}) \subset \mathcal{S}_D^{n+2,a_D}$ and $\tau_{het}(\mathcal{S}_D^{n,a_D}) \subset \mathcal{S}_D^{n,a_D}$.
The final condition is equivalent to $D$ having only classical states. The result follows from Lemma \ref{lemma4.1}. 
\end{proof}

\begin{table}
\centering
{\scriptsize
\begin{tabular}{|c|c|c|c|}
\hline
$s$ & $s^{ij}$ & $a_s-a_{s^{ij}}$ & $\nu_s-\nu_{s^{ij}}$ \\
\hline
$a\,a$ & $A\,A$ & $2$  & $2$  \\
$a\,b$ & $A\,B$ & $1$  & $0$  \\
$a\,A$ & $A\,a$ & $0$  & $0$  \\
$a\,B$ & $A\,b$ & $1$  & $2$  \\
$b\,b$ & $B\,B$ & $0$  & $-2$ \\
$b\,a$ & $B\,A$ & $1$  & $0$  \\
$b\,B$ & $B\,b$ & $0$  & $0$  \\
$b\,A$ & $B\,a$ & $-1$ & $-2$ \\
\hline
\end{tabular}
}
\caption{The first two columns track the change of the labels of the $i\text{-th}$ and $j\text{-th}$ markers from the state $s$ to its transposition $s^{ij}$.}
\label{tab1}
\end{table}

As stated in Section \ref{section2.2}, any two equivalent checkerboard colorable link diagrams can always be related by a finite sequence of generalized Reidemeister moves that preserve the checkerboard colorability (\cite[Theorem\,2.9]{im2}). A direct examination of these moves shows that the property of admitting a classical state is preserved under each move and that, if a checkerboard colorable link diagram does not admit a classical state, then neither does any equivalent checkerboard colorable diagram. We therefore give the following definition

\begin{definition}
A checkerboard colorable link has the \textit{classical-state property} if one (equivalently, every) of its checkeroard colorable diagrams has a classical state.
\end{definition}

To show the following lemma, we recall the notion of \textit{virtual genus} of a virtual link $L$, denoted by $g_v(L)$, and defined as $$g_v(L) = \min \left\lbrace sg(D) | D \text{ is a diagram of } L \right\rbrace.$$ 
As a corollary of \cite[Theorem\,1]{kuperberg}, a virtual link $L$ is classical if and only if $g_v(L)=0$. 

\begin{lemma}
Let $L$ be a checkerboard colorable link. Then, $L$ is classical if and only if $L$ has the classical-state property.
\label{lem4.5}
\end{lemma}

\begin{proof}

Every classical link has the classical-state property. Thus, it suffices to show that a virtual link with the classical-state property is classical. We will do so by showing that the virtual genus is zero. 

Let $d$ be checkerboard colorable link diagram. Since $d$ is alternable, there exists an alternating link diagram $d_{alt}$ obtained from $d$ by finitely many crossing changes. Let $P=(S_d, \widetilde d)$ and $P_{alt}=(S_{d_{alt}}, \widetilde{d_{alt}})$ be the associated abstract link diagrams. As stated in Section \ref{section2.1}, the supporting genera of $P$ and $P_{alt}$ are equal to $g(\Sigma_d)$ and $g(\Sigma_{d_{alt}})$, respectively. Since the underlying graphs $|d|$ and $|d_{alt}|$ are identical, the Carter surfaces $\Sigma_d$ and $\Sigma_{d_{alt}}$ are also identical. Thus $sg(d) = sg(d_{alt})$. Since $g_v(L_d) \leq sg(d) = sg(d_{alt})$, it suffices to show that $sg(d_{alt}) = 0$ to conclude that $L_d$ is classical. Consequently, it is enough to prove the result for alternating virtual links. 

Assume that $L$ is an alternating virtual link with the classical-state property. Let $D$ be an alternating diagram of $L$ with a classical state $s$. The diagram $D$ has the boundary property because it is alternating. Furthermore, the states $s_{\partial}$ and $\overline{s}_{\partial}$ described in Section \ref{section2.2} are obtained from $D$ by resolving all crossings to $\alpha$-smoothing and $\beta$-smoothing, respectively, as stated in \cite[Section\,5]{karimi2018}. We denote by $I_s$ and $E_s$ the sets of internal and external crossings in $s$, respectively. Since all crossings of $D$ have positive Conway sign, $s_{\partial}$ and $\overline{s}_{\partial}$ are obtained from $s$ by changing the smoothing at every crossing in $E_s$ and $I_s$, respectively. Note that changing the smoothing in $s$ at one crossing $c$ splits its single component into two components. If $c'$ is another crossing that does not interlace with $c$ in $s$, then the two parallel arcs in the vicinity of $c'$ lie in the same component after changing the smoothing at $c$. Hence, changing the smoothing at $c'$ splits this component into two, which increases the number of components by $1$. By iteration, we show that if we change the smoothing at $k$ pairwise non-interlacing crossings in $s$, the resulting state has exactly $1+k$ components. Since $s$ is a classical state, the crossings in $I_s$ are pairwise non-interlacing, which is also the case for the crossings in $E_s$. Thus, we have $\# s_{\partial} = 1 + \# I_s$ and $\# \overline{s}_{\partial} = 1 + \# E_s$, where $\# s_{\partial}$ and $\# \overline{s}_{\partial}$ are the numbers of components of the states $s_{\partial}$ and $\overline{s}_{\partial}$, respectively. On the other hand, by \cite[Lemma\,3.18]{karimi2018}, we have $\# s_{\partial} + \# \overline{s}_{\partial} = c(D) + 2 - 2sg(D)$, where $c(D)$ is the number of crossings of $D$. Hence, $(1 + \# I_s) + (1 + \# E_s) = 2 + c(D) = c(D) + 2 - 2sg(D)$, and hence $sg(D) = 0$. Thus $g_v(L) \leq sg(D) = 0$, so $g_v(L) = 0$. Therefore, $L$ is classical.  
\end{proof}

\begin{proof}[Proof of Theorem \ref{th2}]
The result follows directly from the third item of Proposition \ref{prop4.1} and Lemma \ref{lem4.5}.
\end{proof}

\begin{remark}
In general, the determinant and mock-determinant are incomparable. In particular, we provide two checekrboard colorable knots for which the determinant is strictly bigger than the mock-determinant for the first one and vice-versa for the second one. The first knot is the one represented by the top center diagram in Figure \ref{tikz7}, denoted by $D$. This diagram is obtained from the left bottom diagram, denoted by $d$, via a \textit{virtualization} ($\textbf{v}$) at the specified crossing: An operation consisting of first applying a sign change at the specified crossing and then changing the height information. Kauffman showed that the Jones polynomial is invariant under virtualization of crossings \cite[Lemma\,7]{K}. Since the mock-determinant is obtained as the $L_1$-norm of the specialization of the Jones-Kauffman polynomial at the eighth root of unity $e^{i\frac{\pi}{4}}$, the mock-determinant is also invariant under the virtualization of crossings. The diagram $d$ is equivalent to the connected sum of two trefoil knot diagrams $d_1$ and $d_2$ as depicted in the bottom left of Figure \ref{tikz7}. Thus, by \cite[Theorem\,6.4]{boden2022}, we have $\mdet(D)=\mdet(d)=\det(d)= \det(d_1) \times \det(d_2) = 9.$ As we will see in Example \ref{exmpl1}, we have $\det(D)=15$. The second knot is the tabulated knot $4.99$ in Green's table \cite{green}. We obtained in Example \ref{expl3.5} that $\det(4.99) = 3$, and by \cite[Example\,3.8]{abchir2026}, we have $\mdet(4.99) = 5$.
\label{rmrk4.6}
\end{remark}

\begin{remark}
The second item of Proposition \ref{prop4.1} and Lemma \ref{lem4.5} show that classicality of a checkerboard colorable link is sufficient for its determinant and mock-determinant to agree, as expected since these are two extensions of the determinant of classical links. However, classicality is not necessary. Indeed, by using Theorem \ref{th1}, the non-classical checkerboard colorable knot $4.85$ in Green's table \cite{green} satisfies $\det(4.85)=\mdet(4.85)=1$. There are, in fact, infinitely many examples of non-classical knots whose determinants and mock-determinants agree at values different from $1$: Let $K$ be any classical knot with $\det(K) > 1$ and consider the connected sum $K' = 4.85 \# K$. Since checkerboard colorability is preserved under connected sum, the knot $K'$ is checkerboard colorable, while $K'$ is non-classical since it does not have the classical-state property. Since both the determinant and the mock-determinant are multiplicative under the connected sum of checkerboard colorable links (\cite[Theorem\,6.4]{boden2022} and \cite[Proposition\,3.4]{abchir2026}), we have $\mdet(K')= \mdet(4.85) \times \mdet(K) = \mdet(K)$ and $\det(K')= \det(4.85) \times \det(K) = \det(K)$. Since $K$ is classical, $\det(K)=\mdet(K)$, and hence $\mdet(K') = \det(K') > 1$. 
\end{remark}

\section{Quasi-alternating checkerboard colorable links}
\label{sec5}
In the classical setting, the set of quasi-alternating links $\mathcal{Q}$ appeared in the context of link homology as a natural generalization of alternating links. They were defined in \cite{ozsvath} where the authors showed that they are homologically thin for both Khovanov homology and knot Floer homology as alternating links with which they share many properties. On the other hand, it was shown in \cite{ozsvath} that every non-split alternating link is quasi-alternating.

 We introduce the notion of quasi-alternateness for checkerboard colorable links by using the determinant defined by Boden and Karimi \cite{boden2022}. The set of quasi-alternating links is defined recursively as follows:
\begin{definition}
 The set $^{cc}\mathcal{Q}$ of \textit{quasi-alternating links} is the smallest set of checkerboard colorable links satisfying the following properties:
\begin{enumerate}
\item The unknot belongs to $^{cc}\mathcal{Q}$,
\item If $L$ is a link with a diagram $D$ containing a crossing $c$ such that 
\begin{enumerate}
\item the links $L^{c}_{\alpha}$ and $L^{c}_{\beta}$, represented respectively by the smoothings $D^c_{\alpha}$ and $D^c_{\beta}$ of the diagram $D$ at the crossing $c$ as in Figure \ref{fig.1}, are in $^{cc}\mathcal{Q}$ and,
\item $\det(L)=\det(L^{c}_{\alpha})+\det(L^{c}_{\beta}).$
\end{enumerate}
Then $L$ is in $^{cc}\mathcal{Q}$. In this case we say that $c$ is a
\textit{quasi-alternating crossing} of $D$ and that $D$ is \textit{quasi-alternating} at $c$.
\end{enumerate}
\end{definition}

The first immediate observation is that split links are not quasi-alternating due to the fact that they have zero determinant. A second observation is that this definition generalizes the notion of quasi-alternateness from the classical setting to the virtual one in the sense that each quasi-alternating classical link is also in $^{cc}\mathcal{Q}$. On the other hand, there exist non-classical checkerboard colorable links which are quasi-alternating. The following example exhibits one such link. 

\begin{example}
Consider the virtual link diagram $D$ at the top center of Figure \ref{tikz7}. It is easy to check that it is checkerboard colorable and that it represents a knot $K$. We can show that $K$ is non-classical by noting that it does not have the classical-state property.

Smoothing the diagram $D$ at the crossing $c$ results in the two diagrams $D^c_{\alpha}$ and $D^c_{\beta}$ depicted at the top right and top left of Figure \ref{tikz7}, respectively. The links $L^c_{\alpha}$ and $L^c_{\beta}$ are classical since the generalized Reidemeister moves $v_1$ and $v_2$ transform the virtual diagrams $D^c_{\alpha}$ and $D^c_{\beta}$ into classical ones. More precisely, the diagram $D^c_{\beta}$ is equivalent to the non-split alternating link $n(\frac{2}{3} + \frac{2}{3})$, and thus $L^c_{\beta}$ is quasi-alternating. On the other hand, the link diagram $D^c_{\alpha}$ is equivalent to $n(-1 + \frac{2}{3} + \frac{2}{3})$. By \cite[Proposition\,3.1]{champanerkar2015}, the link $L^c_{\alpha}$ is non-split and alternating, and thus it is quasi-alternating. 

We compute the determinants of $D$, $D^c_{\alpha}$, and $D^c_{\beta}$ using Theorem \ref{th1}: $\det(D) = 15$, $\det(D^c_{\alpha}) = 3$, and $\det(D^c_{\beta}) = 12$. Hence, the virtual knot $K$ is non-classical and quasi-alternating.

\begin{figure}
\centering
\includegraphics[width=0.8\linewidth]{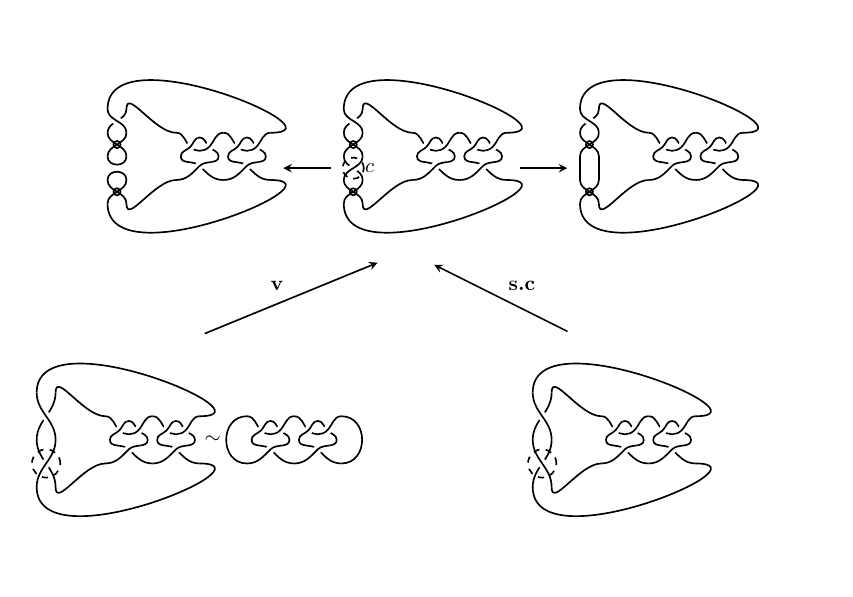}
\caption{}
\label{tikz7}
\end{figure}
\label{exmpl1}
\end{example}

\begin{remark}
Note that the checkerboard colorable diagram $D$ at top center of Figure \ref{tikz7} is obtained by a sign change at a crossing of the classical knot diagram $n(-\frac{1}{2} + \frac{2}{3} + \frac{2}{3})$. More precisely, the sign change is applied to the bottom crossing of the rational tangle $-\frac{1}{2}$ as depicted in Figure \ref{tikz7}, which is a quasi-alternating crossing by \cite[Theorem\,3.8]{abchir2020}. This crossing remains quasi-alternating in $D$ after the sign change as explained in Example \ref{exmpl1}. On the other hand, it is clear that the classical-state property is lost when a sign change is performed at a crossing of a classical link diagram. That is why the virtual link represented by the diagram $D$ is no more classical. This observation can be generalized using the same reasoning as in Example \ref{exmpl1}: If $d$ is a classical link diagram that is quasi-alternating at a crossing $c$, then the checkerboard colorable link diagram $D$ obtained by a sign change on $d$ at $c$ represents a non-classical quasi-alternating link.
\end{remark}

\begin{example}
The checkerboard colorable knot $4.98$ in Green's table \cite{green} is quasi-alternating. Indeed, Figure \ref{qalt} depicts a diagram $D$ of the virtual knot $4.98$ with a specific crossing $c$ such that $D^c_{\beta}$ and $D^c_{\alpha}$ represent the unknot and the alternating virtual link 2-2-3 in Bartholomew's table \cite{bartholomew}, respectively. Furthermore, using Theorem \ref{th1}, we have $\det(4.98)=3$ and $\det(\text{2-2-3})=2$. Hence, $D$ is quasi-alternating at $c$. 
\begin{figure}
\centering
\includegraphics[width=0.8\linewidth]{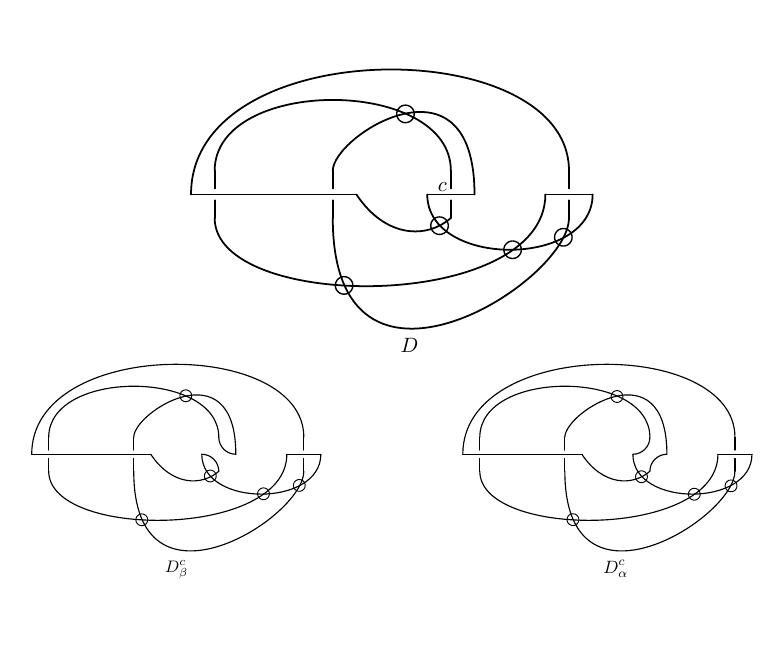}
\caption{}
\label{qalt}
\end{figure}
\label{exmp5.4}
\end{example}

\begin{remark}
By using \ref{eq5}, it is easy to see that non-split alternating virtual links are quasi-alternating. On the other hand, there are non-alternating checkerboard colorable links which are quasi-alternating. For instance, the quasi-alternating knot $4.98$ in the previous example is not alternating.
\end{remark}

\begin{remark}
In \cite[Conjecture\,1.1]{qazaqzeh2013}, the authors conjecture that for every quasi-alternating classical link $L$, the inequality $c(L) \leq \det(L)$ holds, where $c(L)$ the crossing number of $L$. This conjecture was proved for all classical alternating link and infinite families of classical non-alternating quasi-alternating links \cite{qazaqzeh2013,abchir2020}. Based on Example \ref{exmp5.4}, we have $c(4.98) = 4 > \det(4.98) = 3$. This implies that the the aforementioned conjecture cannot be extended to checkerboard colorable quasi-alternating links.
\end{remark} 

Due to the recursive definition, it is almost impossible to show that a given link is not quasi-alternating. An alternative technique to address such a problem is to introduce an obstruction criterion for a checkerboard colorable link to be quasi-alternating. Theorem \ref{th3} provides such an obstruction which will be used subsequently to show the non-quasi-alternateness of some knots in Example \ref{exmpl5.5}. To prove Theorem \ref{th3}, we begin by showing the following lemma.

\begin{lemma}
If $D$ is a checkerboard colorable link diagram with a non-nugatory crossing $c$, then $\mdet(D) \leq \mdet(D^{c}_{\alpha}) + \mdet(D^{c}_{\beta})$.
\label{lem5.6}
\end{lemma}

\begin{proof}
Recall that $\mdet(D) = \|f_D(\xi)\|_1$, where $\xi=e^{i\pi/4}$ and $f_D$ is the Jones-Kauffman polynomial. We show the result in the case where $c$ is of positive Conway sign. The other case is treated similarly. 

If $c$ is a positive crossing, then by \cite[Lemma\,2.5]{abchir2026}, we have 
$$f_D(\xi)=-\xi^{-2}f_{D^{c}_{\alpha}}(\xi)-\xi^{-6e-4}f_{D^{c}_{\beta}}(\xi),$$ which implies

$$
\begin{aligned}
\mdet(D)&= \|f_D(\xi)\|_1 \\
&\leq \|\xi^{-2}f_{D^{c}_{\alpha}}(\xi)\|_1 +\|\xi^{-6e-4}f_{D^{c}_{\beta}}(\xi)\|_1\\
&=\|f_{D^{c}_{\alpha}}(\xi)\|_1+\|f_{D^{c}_{\beta}}(\xi)\|_1 \\
&= \mdet(D^{c}_{\alpha}) + \mdet(D^{c}_{\beta}).
\end{aligned}
$$

If $c$ is a negative crossing, then by \cite[Lemma\,2.5]{abchir2026}, we have 
$$f_L(\xi)=-\xi^{-6e+4}f_{D^{c}_{\beta}}(\xi)-\xi^2f_{D^{c}_{\alpha}}(\xi),$$ and the same argument yields the desired inequality. 
\end{proof}

\begin{proof}[Proof of Theorem \ref{th3}]
Let $L$ be a quasi-alternating checkerboard colorable link. We use induction on $\det(L)$. If $\det(L) = 1$, then $L$ is the unknot, and thus $\mdet(L) = 1$. This proves the induction basis. Suppose the result holds for every quasi-alternating checkerboard colorable link with a determinant less that $\det(L)$. For a checkerboard colorable diagram $D$ of $L$ that is quasi-alternating at a crossing $c$, the induction hypothesis implies that $\mdet(D^{c}_{\alpha}) \leq \det(D^{c}_{\alpha})$ and $\mdet(D^{c}_{\beta}) \leq \det(D^{c}_{\beta})$. By Lemma \ref{lem5.6}, we get $\mdet(D) \leq \det(D^{c}_{\alpha}) + \det(D^{c}_{\beta}) = \det(D)$. The result then follows by induction.
\end{proof}

\begin{example}

The checkerboard colorable knots $3.5$ and $4.99$ in Green's table \cite{green} were discussed in \cite[Example\,4.2]{abchir2026} where it was shown that $\mdet(3.5)=3$ and $\mdet(4.99)=5$. On the other hand, using Theorem \ref{th1}, we find that $ \det(3.5) = 1$ and $\det(4.99) = 3$. Thus, by Theorem \ref{th3}, both $3.5$ and $4.99$ are not quasi-alternating.
\label{exmpl5.5}
\end{example} 

\begin{remark}
In \cite[Example\,4.2]{abchir2026}, it was shown that the checkerboard colorable knots $3.5$ and $4.99$ belong to the set of mock-alternating links introduced in \cite{abchir2026}, which is defined analogously to the set of quasi-alternating checkerboard colorable links, where the mock-determinant is used instead of the determinant. Since every alternating link is quasi-alternating, the observation in \cite[Remark\,4.7]{abchir2026} that not every alternating link is mock-alternating shows that quasi-alternateness does not imply mock-alternateness. On the other hand, as shown in the previous example, the checkerboard colorable knots $3.5$ and $4.99$ are not quasi-alternating. Thus, mock-alternateness does not imply quasi-alternateness either.
\end{remark}

\subsection{Twisting checkerboard colorable virtual links}

In this section, we use the rational extension technique to give rise to an infinite family of quasi-alternating links from a given quasi-alternating one. This technique was described  in \cite{kofman2009} for the first time in the classical setting using rational tangles and later generalized in \cite{qazaqzeh2015} to a product of rational tangles and to non-algebraic alternating tangles in \cite{chbili2020}. The basic idea of this technique lies on replacing a quasi-alternating crossing $c$ by some tangle. 

\begin{theorem}
If $D$ is a quasi-alternating diagram at some crossing $c$, then the diagram $D^{c\leftarrow t}$ obtained from $D$ by replacing the crossing $c$ by a rational tangle $t$ that extends $c$ is quasi-alternating diagram at any crossing of the tangle $t$.
\label{new}
\end{theorem}

\begin{proof}
By changing the checkerboard coloring, taking the mirror image, or both, if necessary, we may assume that the crossing $c$ in the diagram $D$ satisfies $\varepsilon(c) = +1$ and has a negative Conway sign. We prove the result by induction on the length of the continued fraction of the given rational tangle. If the length is one, then we have an integer tangle in such a case. Now we prove this case by induction on the number of crossings of this integer tangle. Both such inductions depend on proving the result for the diagrams $d:=D^{c \leftarrow \frac{1}{2}}$ and $d' :=D^{c \leftarrow 2}$ where the crossing $c$ has been replaced, respectively, by a vertical or horizontal integer tangle of two crossings that extends $c$. Now we consider the first diagram $d$. If we smooth the second crossing, which by abuse of notation we also denote by $c$ in the new diagram $d$, we obtain the links $d^c_{\alpha}$ and $d^c_{\beta}$ that represent the diagrams $D$ and $D^c_{\beta}$, respectively. It is clear that these are checkerboard colorable quasi-alternating diagrams from the fact that $D$ is checkerboard colorable quasi-alternating at the crossing $c$. To prove the result, it is left to prove that $\det(d) = \det(d^c_{\alpha})  + \det(d^c_{\beta})$. We have $\det(d^c_{\alpha})  + \det(d^c_{\beta}) = \det(D) + \det(D^c_{\beta}) = \det(D^c_{\alpha}) + 2\det(D^c_{\beta}).$ By Proposition \ref{prop3.11}, we have $\det(d) = \det(d^c_{\alpha})  + \det(d^c_{\beta})$. Thus, we conclude that $d$ is quasi-alternating. The same can similarly be proven for the diagram $d'$. Now we assume that the result holds for $D^{c \leftarrow \frac{1}{n-1}}$ and that it is quasi-alternating at the last crossing of the vertical tangle $\frac{1}{n-1}$. Now we incorporate the previous argument to show the result for $D^{c \leftarrow \frac{1}{n}}$. Similarly, we obtain that the result holds for $D^{c \leftarrow n}$ also. Suppose now that the result holds for any rational tangle extending $c$ whose fraction has length $n$.Note that $[a_{1},a_{2}, \ldots, a_{n}+1]=[a_{1},a_{2}, \ldots, a_{n},1]$ and that $[a_{1},a_{2}, \ldots, a_{n+1}]$ is obtained from $[a_{1},a_{2}, \ldots, a_{n},1]$ by extending the last crossing with the integer tangle $a_{n+1}$ which can be vertical or horizontal depending the the parity of $n$ (see Figure \ref{rt}). At the end, by using the first induction on the number of crossings in the integer tangle $a_{n+1}$, we obtain that the result holds for any rational tangle that extends the crossing $c$.
\end{proof}

\subsection{Dealternator extensions of almost alternating links}

In this paragraph, we provide a quasi-alternating link by a rational extension at the dealternator of an almost alternating diagram representing a quasi-alternating link. We begin by two observations which easily follow from Proposition \ref{prop2}.

\begin{corollary}
Let $D$ be an almost alternating virtual link diagram with dealternator $c$. Then, we have the following
\begin{enumerate}
\item[(i)] $D$ is never quasi-alternating at $c$.
\item[(ii)] If the link represented by $D$ is quasi-alternating, then $\det(D^c_{\alpha}) \neq \det(D^c_{\beta})$.
\end{enumerate}
\end{corollary}

Theorem \ref{new} states that the twisting technique preserves quasi-alternateness when applied at a quasi-alternating crossing. Despite the fact that dealternator is never a quasi-alternating crossing, we prove that if an almost alternating virtual link is quasi-alternating, then the twisting technique, when applied at the dealternator, also preserves quasi-alternateness as stated in the following result.

\begin{theorem}
  Let $D$ be an almost alternating diagram representing a quasi-alternating link. If $c$ is the dealternator of the diagram $D$, then we have the following:
  \begin{enumerate}
  \item If $\det ( D^{c}_{\alpha} ) > \det ( D^{c}_{\beta} )$, then the diagram $D^{c \leftarrow -\frac{a}{b}}$ is quasi-alternating if $\dfrac{a}{b} < 1$.
\item If $\det ( D^{c}_{\alpha} ) < \det ( D^{c}_{\beta} )$, then the diagram $D^{c \leftarrow -\frac{a}{b}}$ is quasi-alternating if $\dfrac{a}{b} > 1$.
\end{enumerate}
\label{th5.5}
\end{theorem}

\begin{proof}
As explained in Remark \ref{rmrk2}, we can write $D$ in the form $n(-1 + T)$ where the dealternator $c$ is identified with the rational tangle $-1$ and $T$ is an alternating virtual tangle diagram of type $1$. We simply show the first result and the second one can be obtained by a similar argument. 

Assume that $\det ( D^{c}_{\alpha} ) > \det ( D^{c}_{\beta} )$, or equivalently, $\det ( n(T) ) > \det ( d(T) )$. We denote the link diagram $D^{c \leftarrow -\frac{1}{2}}$ by $d$. By abuse of notation, we denote by $c$ the top crossing of the vertical tangle $-\frac{1}{2}$. We first show that $d$ is quasi-alternating at $c$. It is clear that $d^{c}_{\alpha}$ is equivalent to the alternating diagram $D^{c}_{\alpha} = n(T)$, and that $d^{c}_{\beta}$ is equivalent to $D$. It remains to show that $\det(d) = \det(n(T)) + \det(D)$ which is easily deduced from the assumption $\det ( n(T) ) > \det ( d(T) )$ and the second item of Proposition \ref{prop3.10}. Hence, $d$ is quasi-alternating at $c$. Now, it suffices to observe that $D^{c \leftarrow -\frac{a}{b}} = d^{c \leftarrow -\frac{a}{b-a}}$, and the result follows from Theorem \ref{new}. 
\end{proof}

\subsection*{Acknowledment}
The authors acknowledge that AI has been used for the sole purpose of typesetting support, grammatical correction, and paraphrasing  to enhance readability and improve the manuscript's language. After this usage, the authors reviewed, edited and approved the final content, and took full responsibility for the scientific integrity and accuracy of the text.

\end{document}